\documentclass[letterpaper, final]{IEEEtran} 
\IEEEoverridecommandlockouts 
\usepackage{mathtools} 
\usepackage{hyperref,graphicx} 
\usepackage{times} 
\usepackage{amsmath} 
\usepackage{amssymb}  
\usepackage{multirow}
\usepackage{color}
\usepackage{algorithm, algorithmicx, algpseudocode}

\newenvironment{proof-of}[1]{\quad{\it Proof of #1.\,}}{\hfill \IEEEQED \par}

\newenvironment{pf}{\quad{\it Proof.\,}}{\hfill \IEEEQED \par}
\newtheorem{thm}{Theorem}
\newtheorem{cor}{Corollary}

\newtheorem{defn}{Definition}
\newtheorem{lem}{Lemma}

\newtheorem{prop}{Proposition}

\newtheorem{rem}{Remark}

\newcommand{\R}{{\mathbb R}}
\newcommand{\Rnn}{{\mathbb R}_{\ge 0}}
\newcommand{\Rp}{{\mathbb R}_{> 0}}

\newcommand{\cA}{{\mathcal A}}
\newcommand{\cB}{{\mathcal B}}

\newcommand{\cD}{{\mathcal D}}

\newcommand{\cF}{{\mathcal F}}
\newcommand{\cG}{{\mathcal G}}

\newcommand{\cM}{{\mathcal M}}

\newcommand{\cP}{{\mathcal P}}

\newcommand{\cS}{{\mathcal S}}
\newcommand{\cT}{{\mathcal T}}
\newcommand{\cU}{{\mathcal U}}
\newcommand{\cV}{{\mathcal V}}
\newcommand{\cW}{{\mathcal W}}
\newcommand{\cX}{{\mathcal X}}

\newcommand{\diag}{{\mathrm{diag}}}

\newcommand{\cl}{\mathrm{cl}}

\newcommand{\bfone}{\mathbf 1}
\newcommand{\inter}{\mathrm{int}}

\title{Shaping Pulses to Control Bistable Biological Systems.}
\author{ Aivar Sootla,  Diego Oyarz\'{u}n,  David Angeli and Guy-Bart Stan
\thanks{
Dr Sootla is with Institut Montefiore, University of Li\`{e}ge, B-4000, Belgium {\tt asootla@ulg.ac.be}. Dr Oyarz\'{u}n, Dr Angeli, Dr Stan {\tt (d.oyarzun, d.angeli, g.stan)@imperial.ac.uk} are with Departments of Mathematics, Electrical and Electronic Engineering, and Bioengineering, respectively, Imperial College London, SW72AZ, London, UK. Dr Angeli is also with the Department of Information Engineering, University of Florence, Italy.

Part of this work was performed, while Dr Sootla was a post-doctoral research associate at ICL. Dr Sootla and Dr Stan acknowledge support by the EPSRC grants EP/J014214/1, EP/G036004/1. Dr Stan is additionally supported by the EPRSC grant EP/M002187/1. Dr Sootla is now supported by an F.R.S.--FNRS Fellowship. Dr Oyarz\'{u}n is supported by a Junior Research Fellowship from ICL. The authors would like to thank Dr Alexandre Mauroy for valuable suggestions and discussions. This technical report contains the material from publications~\cite{sootla2015pulsesacc} and~\cite{sootla2015pulsesaut}.}}
\begin{document}
\maketitle
\begin{abstract}	
In this paper we study how to shape temporal pulses to switch a bistable system between its stable steady states. Our motivation for pulse-based control comes from applications in synthetic biology, where it is generally difficult to implement real-time feedback control systems due to technical limitations in sensors and actuators. We show that for monotone bistable systems, the estimation of the set of all pulses that switch the system reduces to the computation of one non-increasing curve. We provide an efficient algorithm to compute this curve and illustrate the results with a genetic bistable system commonly used in synthetic biology. We also extend these results to models with parametric uncertainty and provide a number of examples and counterexamples that demonstrate the power and limitations of the current theory. In order to show the full potential of the framework, we consider the problem of inducing oscillations in a monotone biochemical system using a combination of temporal pulses and event-based control. Our results provide an insight into the dynamics of bistable systems under external inputs and open up numerous directions for future investigation. 
\end{abstract}
\section{Introduction}
In this paper we investigate how to switch a bistable system between its two stable steady states using external input signals. Our main motivation for this problem comes from synthetic biology, which aims to engineer and control biological functions in living cells~\cite{brophy2014principles}. Most of current research in synthetic biology focusses on building biomolecular circuits inside cells through genetic engineering. Such circuits can control cellular functions and implement new ones, including cellular logic gates, cell-to-cell communication and light-responsive behaviours. These systems have enormous potential in diverse applications such as metabolic engineering, bioremediation, and even the energy sector~\cite{Purnick:2009}. 

Several recent works~\cite{milias2011silico,Menolascina:2011,uhlendorf2012long} have showcased how cells can be controlled externally via computer-based feedback  and actuators such as chemical inducers or light stimuli \cite{Mettetal08,Levskaya09}. An important challenge in these approaches is the need for real-time measurements, which tend to be costly and difficult to implement with current technologies. In addition, because of technical limitations and the inherent nonlinearity of biochemical interactions, actuators are severely constrained in the type of input signals they can produce. As a consequence, the input signals generated by traditional feedback controllers (e.g. PID or model predictive control) may be hard to implement without a significant decrease in control performance.

In this paper we show how to switch a bistable system without the need for output measurements. We propose an open-loop control strategy based on a temporal pulse of suitable magnitude $\mu$ and duration $\tau$: 
\begin{align}
	u(t) &= \mu h(t, \tau),\qquad  h(t,\tau)= \begin{cases}
	1 & 0\leq t \leq \tau,\\
	0 & t> \tau.
	\end{cases}\label{eq:pulse}
\end{align} 
Our goal is to characterise the set of all pairs $(\mu,\tau)$ that can switch the system between the stable steady states and the set of all pairs $(\mu,\tau)$ that cannot. We call these sets {\it the switching sets} and a boundary between these sets the \emph{switching separatrix}. The pairs $(\mu,\tau)$ close to the switching separatix are especially important in synthetic biology applications, as a large $\mu$ or a large $\tau$ can trigger toxic effects that slow down cell growth or cause cell death. 

In a previous paper~\cite{sootla2015pulsesacc}, we showed that for monotone systems the switching separatrix is a monotone curve. This result was therein extended to a class of non-monotone systems whose vector fields can be bounded by vector fields of monotone systems. This idea ultimately leads to robustness guarantees under parametric uncertainty. These results are in the spirit of~\cite{gennat2008computing,ramdani2010computing,ramdani2009hybrid}, where the authors considered the problem of computing reachability sets of a monotone system. Some parallels can be also drawn with~\cite{meyer2013controllability,chisci2006asymptotic}, where feedback controllers for monotone systems were proposed. 

\emph{Contributions.} In the present paper we provide the first complete proof of our preliminary results in~\cite{sootla2015pulsesacc} and extend them in several directions. We formulate necessary and sufficient conditions for the existence of the monotone switching separatrix for non-monotone systems. Although it is generally hard to use this result to establish monotonicity of the switching separatrix, we use it to prove the converse. For example, we show that for a bistable Lorenz system the switching separatrix is not monotone. We then generalise the main result of~\cite{sootla2015pulsesacc} by providing conditions for the switching separatrix to be a graph of a function. We also discuss the relation between bifurcations and the mechanism of pulse-based switching, which provides additional insights into the switching problem. We use this intuition to show and then explain the failure of pulse-based control on an HIV viral load control problem~\cite{adams2004dynamic}. We proceed by providing a numerical algorithm to compute the switching separatrices for monotone systems. The algorithm can be efficiently distributed among several computational units and does not explicitly use the vector field of the model. We evaluate the theory and computational tools on the bistable LacI-TetR system, which is commonly referred to as a \emph{genetic toggle switch}~\cite{Gardner00}. Genetic toggle switches are widely used in synthetic biology to trigger cellular functions in response to extracellular signals~\cite{brophy2014principles,khalil2010synthetic}.

We complement our theoretical findings with several observations that illustrate limitations of the current theory and highlight the need for deeper investigations of bistable systems. For example, we show that for a toxin-antitoxin system~\cite{cataudella2013conditional}, the switching separatrix appears to be monotone, even though the system does not appear to be monotone. Finally, in order to demonstrate the full potential of pulse-based control, we consider the problem of inducing an oscillatory behaviour in a generalised repressilator system~\cite{Strelkowa10}. 

\emph{Organisation.} The rest of the paper is organised as follows. In Section~\ref{s:prem} we cover the basics of monotone systems theory, formulate the problem in Subsection~\ref{ss:pf}, and provide an intuition into the mechanism of pulse-based switching for monotone systems in Subsection~\ref{ss:mech}. We also provide some motivational examples for the development of our theoretical results.
In Section~\ref{s:theory} we formulate the theoretical results and in Section~\ref{s:comp} we present the computational algorithm, which we evaluate in Section~\ref{s:laci-tetr} on the LacI-TetR system. In Section~\ref{s:examples}, we provide counterexamples and an application of inducing oscillations in a generalised repressilator system. All the proofs are found in the Appendix.

{\it Notation.} Let $\|\cdot\|_2$ stand for the Euclidean norm in $\R^n$, $Y^{\ast}$ stand for a topological dual to $Y$, $X\backslash Y$ stand for the relative complement of $X$ in $Y$, $\inter(Y)$ stand for the interior of the set $Y$, and $\cl(Y)$ for its closure. 
\section{Preliminaries\label{s:prem}} 
Consider a single input control system
\begin{equation}
\label{eq:sys}
\dot x = f(x,u),\quad x(0) = x_0,
\end{equation} 
where $f: \cD\times \cU\rightarrow \R^n$, $u:\R_{\ge 0}\rightarrow \cU$, $\cD\subset\R^n$, $\cU\subset\R$ and $u(\cdot)$ belongs to the space $\cU_{\infty}$ of Lebesgue measurable functions with values from $\cU$. We say that the system is \emph{unforced}, if $u=0$. We define \emph{the flow} map $\phi_f: \R \times \cD \times \cU_{\infty}\rightarrow \R^n$, where $\phi_f(t, x_0, u)$ is a solution to the system~\eqref{eq:sys} with an initial condition $x_0$ and a control signal $u$. We consider the control signals in the shape of \emph{a pulse}, that is signals  defined in~\eqref{eq:pulse} with the set of admissible $\mu$ and $\tau$ denoted as $\cS = \{\mu, \tau \in \Rnn\}$. 

In order to avoid confusion, we reserve the notation $f(x,u)$ for the vector field of non-monotone systems, while systems
\begin{align}
\label{sys:low}	\dot x &= g(x,u), &x(0) &= x_0,\\
\label{sys:up} \dot x &= r(x,u), &x(0) &= x_0,
\end{align}
denote so-called \emph{monotone systems} throughout the paper. In short, monotone systems are those which preserve a partial order relation in initial conditions and input signals. A relation $\sim$ is called a {\it partial order} if it is reflexive ($x\sim x$), transitive ($x\sim y$, $y\sim z$ implies $x\sim z$), and antisymmetric ($x\sim y$, $y\sim x$ implies $x = y$). We define a partial order $\succeq_x$ through a cone $K\in\R^n$ as follows: $x\succeq_x y$ if and only if $ x - y \in K$. We write $x\not \succeq_x y$, if the relation  $x \succeq_x y$ does not hold. We will also write $x\succ_x y$ if $x\succeq_x y$ and $x\ne y$, and $x\gg_x y$ if $x- y \in \inter(K)$. Similarly we can define a partial order on the space of signals $u\in \cU_{\infty}$: $u\succeq_u v$ if $u(t) - v(t) \in K$ for all $t\ge 0$. We write $u\succ_u v$, if $u\succeq_u v$ and $u(t) \ne v(t)$ for all $t\ge 0$. Partial orders also induce some geometric properties on sets. A set $M$ is called \emph{p-convex} if for every $x$, $y$ in $M$ such that $x\succeq_x y$, and every $\lambda\in(0,1)$ we have that $\lambda x+ (1-\lambda) y\in M$.

\begin{defn}
 The system~\eqref{sys:low} is called \emph{monotone} on $\cD_M\times \cU_{\infty}$ with respect to the partial orders $\succeq_x$, $\succeq_u$, if for all $x, y\in\cD_M$ and $u, v\in\cU_{\infty}$ such that $x\succeq_x y$ and $u \succeq_u v$, we have $\phi_g(t, x, u) \succeq_x \phi_g(t,y, v)$ for all $t\ge 0$. If additionally, $x\succ_x y$, or $u \succ_x v$ implies that $\phi_g(t, x, u) \gg_x \phi_g(t,y, v)$ for all $t>0$, then the system is called \emph{strongly monotone}.
\end{defn}

In general, it is hard to establish monotonicity of a system with respect to an order other than an order induced by an orthant (e.g., positive orthant $\Rnn^n$). Hence throughout the paper, by a monotone system we actually mean \emph{a monotone system with respect to a partial order induced by an orthant}. A certificate for monotonicity with respect to an orthant is referred to as Kamke-M\"uller conditions~\cite{angeli2003monotone}.
\begin{prop}[\cite{angeli2003monotone}]\label{prop:kamke}
Consider the system~\eqref{sys:low}, where $g$ is differentiable in $x$ and $u$ and let the sets $\cD_M$, $\cU$ be p-convex. Let the partial orders $\succeq_x$, $\succeq_u$ be induced by $P_x \Rnn^n$, $P_u \Rnn^m$, respectively, where $P_x = \diag((-1)^{\varepsilon_1}, \dots, (-1)^{\varepsilon_n})$, $P_u = \diag((-1)^{\delta_1}, \dots, (-1)^{\delta_m})$ for some $\varepsilon_i$, $\delta_i$ in $\{0 ,1\}$. Then 
\begin{gather*}
(-1)^{\varepsilon_i + \varepsilon_j}\frac{\partial g_i}{\partial x_j}\ge 0,\quad\forall~i\ne j,\quad(x,u)\in\cl(\cD_M)\times\cU\\
(-1)^{\varepsilon_i+\delta_j}\frac{\partial g_i}{\partial u_j}\ge 0,\quad\forall~i, j,\quad(x,u)\in\cD_M\times\cU
\end{gather*}
if and only if the system~\eqref{sys:low} is monotone on $\cD_M\times \cU_{\infty}$ with respect to $\succeq_x$, $\succeq_u$
\end{prop}
If we consider the orthants $\Rnn^n$, $\Rnn^m$, then the conditions above are equivalent to checking if for all $x \preceq_x y$ such that $x_i = y_i$ for some $i$, and all $u\preceq_u v$ we have $g_i(x,u) \le g_i(y,v)$.

\subsection{Problem Formulation\label{ss:pf}}
We confine the class of considered control systems by making the following {\bf assumptions:}
  \begin{enumerate}
  \item[{\bf A1.}] Let $f(x,u)$ in~\eqref{eq:sys} be continuous in $(x,u)$ on $\cD_f\times \cU$. Moreover, for each compact sets $C_1\subset \cD_f$ and $C_2 \subset \cU$, let there exist a constant $k$ such that $\|f(\xi, u) - f(\zeta, u)\|_2 \le k \|\xi - \zeta\|_2$ for all $\xi,\zeta \in C_1$ and $u\in C_2$.
  \item[{\bf A2.}] Let the unforced system~\eqref{eq:sys} have two stable steady states in $\cD_f$, denoted as $s^0_f$ and $s^1_f$,
  \item[{\bf A3.}] Let $\cD_f=\cl(\cA(s_{f}^{0})\cup \cA(s_{f}^{1}))$, where $\cA(s^i_f)$ stands for the domain of attraction of the steady state $s^i_f$ for $i=0,1$ of the unforced system~\eqref{eq:sys}, 
  \item[{\bf A4.}] For any $u = \mu h(\cdot, \tau)$ with finite $\mu$ and $\tau$ let $\phi_f(t, s_f^0, u)$ belong to $\inter(\cD_f)$. Moreover, let the sets 
  \begin{align*}
\cS^+_f &= \{\mu, \tau>0\Bigl| \lim_{t\to\infty} \phi_{f}(t;s_{f}^{0},\mu h(\cdot,\tau)) = s_{f}^{1}\}\\
\cS^-_f &= \{\mu, \tau>0\Bigl| \lim_{t\to\infty} \phi_{f}(t;s_{f}^{0},\mu h(\cdot,\tau)) = s_{f}^{0}\}
\end{align*} 
be non-empty.
\end{enumerate}

Assumption~A1 guarantees existence, uniqueness and continuity of solutions to~\eqref{eq:sys}, while Assumptions~A2--A3 define a bistable system on a set $\cD_f$ controlled by pulses. In Assumption~A4 we define the \emph{switching sets}: the set $\cS^+_f$, which contains all $(\mu,\tau)$ pairs that switch the system, and the set $\cS^-_f$, which contains all pairs that do not. The boundary between these sets is called the \emph{switching separatrix}. In the rest of the paper, we focus on the {\bf control problem} of estimating the switching sets. 

\subsection{Mechanism of Pulse-Based Switching\label{ss:mech}}
The general problem of switching a bistable system with external inputs is amenable to an optimal control formulation. However, in applications such as synthetic biology, optimal control solutions can be very hard to implement due to technical limitations in actuators and output measurements. Additionally, the solution of this optimal control problem may be technically challenging. Hence applying open-loop pulses can be a reasonable solution, if we can guarantee some form of robustness. As we shall see later, our results show that for monotone systems, pulse-based switching is computationally tractable and robust towards parameter variations. 

Before presenting our main results, we first provide an intuitive link between monotonicity and the ability to switch a system with temporal pulses. If we consider constant inputs $u=\mu$ and regard $\mu$ as a bifurcation parameter, we have the following result.

\begin{prop}\label{prop:bif} Let the system~\eqref{sys:low} satisfy Assumptions~A1--A4 and be monotone on $\cD_g\times\cU_{\infty}$ with respect to $\Rnn^n$, $\Rnn$. Let $\mu_{\rm min}$ be such that all pairs $(\mu,\tau)\in\cS_g^-$ for $0<\mu<\mu_{\rm min}$, and any finite positive $\tau$. Let also $\xi(\mu) = \lim\limits_{t\rightarrow\infty}\phi_g(t, s^0_g, \mu)$ and $\eta(\mu)= \lim\limits_{t\rightarrow\infty}\phi_g(t, s^1_g, \mu)$. Then
\begin{enumerate}
\item If $\mu\le \lambda<\mu_{\rm min}$ then $\xi(\mu)\preceq_x \xi(\lambda)$, $\eta(\mu)\preceq_x \eta(\lambda)$ 
\item If $0<\mu<\mu_{\rm min}$ then $\xi(\mu)\in\cA(s^0_g)$ and $\xi(\mu) \prec_x \eta(\mu)$.
\item The function $\xi(\mu)$ is discontinuous at $\mu_{\rm min}$.
\end{enumerate}
\end{prop}
\begin{figure}[t]\centering
  \includegraphics[height = 0.34\columnwidth]{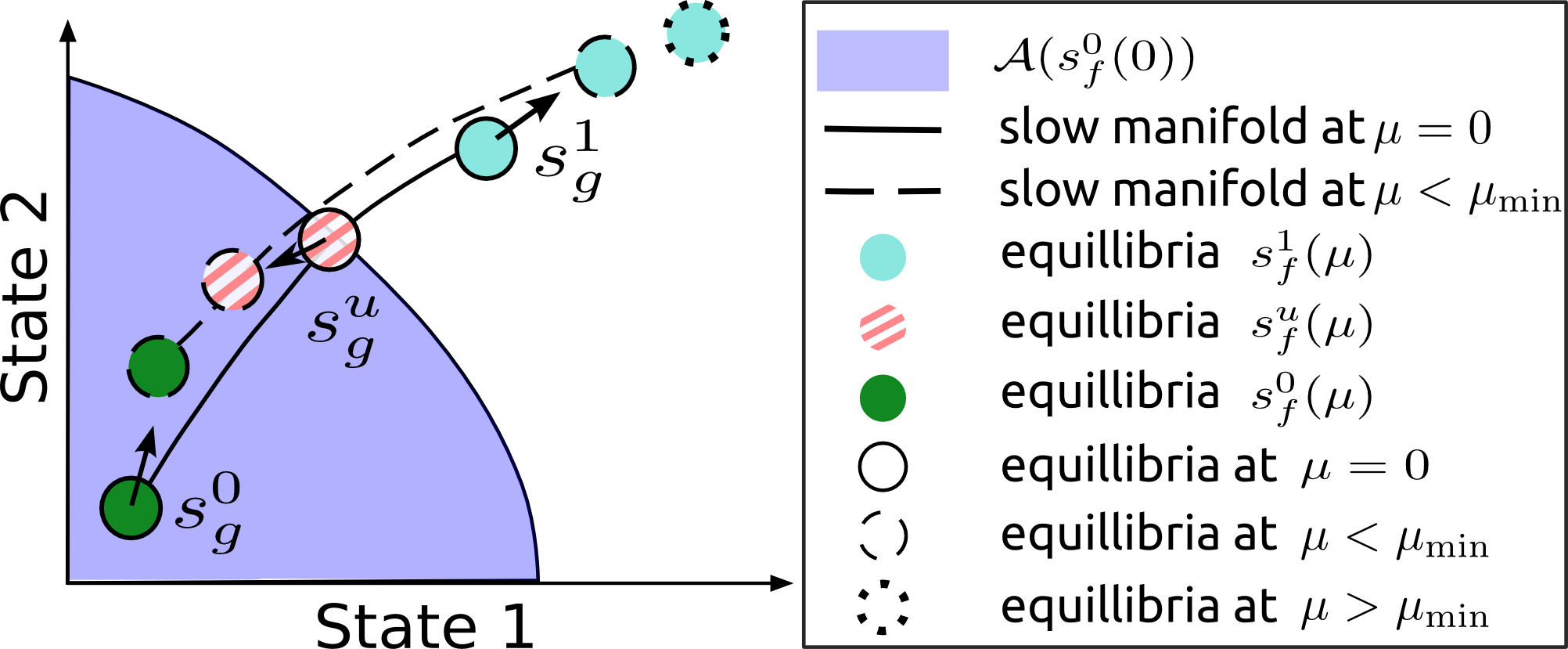}
\caption{A schematic depiction of the evolution of the stable nodes $s_g^0(\mu)$, $s_g^1(\mu)$ and the saddle $s_g^u(\mu)$ with respect to $\mu$ in the genetic toggle switch system. By slow manifold we mean a manifold connecting stable equillibria and a saddle. The arrows show the direction of the equillibria movements with increasing $\mu$. At $\mu_{\rm min}$ the equilibria $s_g^0(\mu)$ and $s_g^u(\mu)$ collide resulting in a saddle-node bifurcation preserving only $s_g^1(\mu)$. } \label{fig:bif}
\end{figure}

The proof of the proposition is in the Appendix. In many applications, the functions $\xi(\mu)$, $\eta(\mu)$ are simply evolutions of the steady states $s^0_g$, $s^1_g$ with respect to the parameter $\mu$, respectively. Hence, statement (1) of Proposition~\ref{prop:bif} shows how the steady states move with respect to changes in $\mu$. Statement (2) ensures that there are at least two distinct asymptotically stable equilibria for $\mu<\mu_{\rm min}$. Finally, statement (3) indicates that the system undergoes a bifurcation for $\mu=\mu_{\rm min}$. The particular type of the bifurcation will depend on a specific model. Next we investigate further aspects of this result with some examples of monotone and non-monotone bistable systems. 

\emph{Example 1: LacI-TetR Switch.} The genetic toggle switch is composed of two mutually repressive genes \emph{LacI} and \emph{TetR} and was a pioneering genetic system for synthetic biology~\cite{Gardner00}.  We consider its control-affine model, which is consistent with a toggle switch actuated by light induction~\cite{Levskaya09}: 
  \begin{equation}
    \label{eq:ts-par}
    \begin{aligned}
      \dot x_1 &= \frac{p_1}{1 + (x_2/p_2)^{p_3}} + p_4 - p_5 x_1 + u, \\
      \dot x_2 &= \frac{p_6}{1 + (x_1/p_7)^{p_8}} + p_9 - p_{10} x_2, 
    \end{aligned}
  \end{equation}
  where the parameters have the following values 
  \begin{equation}\label{eq:ts-par-val}
  \begin{aligned}
   p_1 &= 40, &p_2&=1, &p_3&=4, &p_4&= 0.05, &&p_5   =1, \\
   p_6 &= 30, &p_7&=1, &p_8&=4, &p_9&= 0.1,  &&p_{10}=1.
  \end{aligned}
  \end{equation}
In the model \eqref{eq:ts-par}, $x_{i}$ represents the concentration of each protein, whose mutual repression is modelled via a rational function. The parameters $p_2$ and $p_7$ represent the repression thresholds, whereas $p_4$ and $p_9$ model the basal synthesis rate of each protein. The parameters $p_5$ and $p_{10}$ are the degradation rate constants,  $p_3$, $p_8$ are called Hill (or cooperativity) parameters, and $p_1$, $p_6$ describe the strength of mutual repression. By means of Proposition~\ref{prop:kamke} we can readily check that the model is monotone on $\Rnn^2\times\Rnn$ for all nonnegative parameter values. It can be verified by direct computation that the system satisfies Assumptions~A1--A4 with $\cD_f = \Rnn^2$. It can be also shown that the unforced system is strongly monotone in $\inter(\Rnn^2)$ using the results in~\cite{smith2008monotone}. 

With the chosen parameter values, we numerically found a bifurcation to occur at $\mu_{\rm min} \approx 1.4077$. For $\mu<\mu_{\rm min}$ the system has two stable nodes and a saddle. We observe that $\xi(\mu)=\eta(\mu)$ for all $\mu >\mu_{\rm min}$, and therefore we conclude that the system undergoes a saddle-node bifurcation, as illustrated in Figure~\ref{fig:bif}.

\begin{figure}[t]\centering
  \includegraphics[width = 0.7\columnwidth]{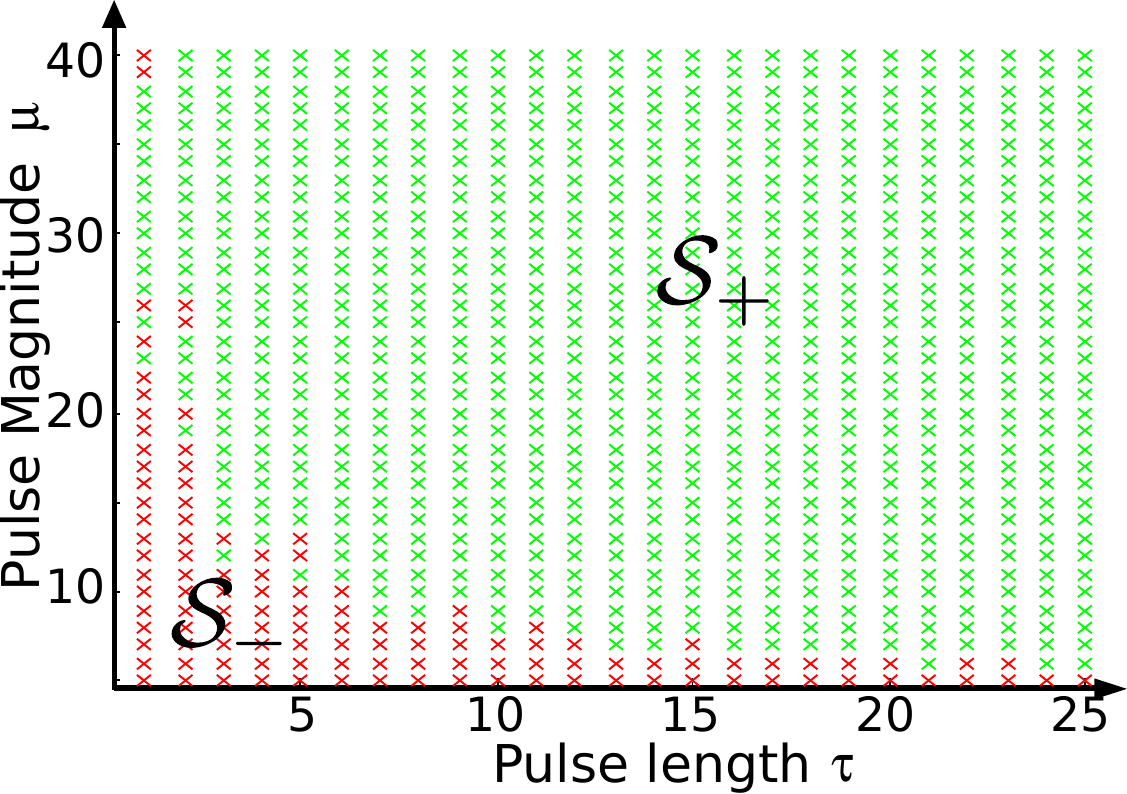}
\caption{Switching sets for the Lorenz system. We simulated the Lorenz system for $(\mu,\tau)$ pairs taken from a mesh grid. The green and red crosses correspond to the pairs that switched or not switched the system, respectively. } \label{fig:lorenz-ss}
\end{figure}

\emph{Example 2: Lorenz system.} Consider a system 
\begin{align*}
\dot x_1 &= \sigma (x_2 - x_1)  + u\\
\dot x_2 &= x_1 (\rho - x_3) - x_2 + u\\
\dot x_3 &= x_1 x_2 -\beta x_3
\end{align*}
with parameters $\sigma = 10$, $\rho = 21$, $\beta = 8/3$, which is non-monotone and bistable with two stable foci. Numerical computation of the  sets $\cS^-$ and $\cS^+$ in Figure~\ref{fig:lorenz-ss} suggests that the switching separatrix is not monotone. We will revisit this conclusion in the next section using our theoretical results.

\emph{Example 3: HIV viral load control problem.} In~\cite{adams2004dynamic} the authors considered the problem of switching from a ``non-healthy'' ($s^0$) to a ``healthy'' ($s^1$) steady state by means of two control inputs ($u_{1}$ and $u_{2}$) that model different drug therapies. Due to space limitations we refer the reader to~\cite{adams2004dynamic} for a description of the model. It can be shown that both steady states are stable foci and that the model is non-monotone. Although the system can be switched with non-pulse control signals~\cite{adams2004dynamic}, using extensive simulations we were unable to find a combination of pulses in $u_{1}$ and $u_{2}$ switching the system. 

As in the case of monotone bistable system, we found a bifurcation with respect to constant control signals $u_1 = \mu_1$ and $u_2= \mu_2$. More specifically, we fixed $\mu_2 = 0.4$, and numerically found a bifurcation at $\mu_1 \approx 0.7059$. The major difference between this case and the monotone system case (Example 1) is that the steady state $s^1(0.7059, 0.4)$ lies the domain of attraction of $s^0(0,0)$. Hence if we stop applying the constant control signal we regress back to the initial point $s^0(0,0)$. Furthermore, with increasing $\mu_1$ the steady state $s^1(\mu_1, 0.4)$ is moving towards the origin, which also lies in the domain of attraction of $s^0(0,0)$. This makes pulse-based switching very difficult, if not impossible. 

\section{Theoretical Results\label{s:theory}}
In~\cite{sootla2015pulsesacc} we showed that the switching separatrix of a monotone bistable system $\dot x = g(x,u)$ is non-increasing. Here we present a generalisation of this result by formulating necessary and sufficient conditions for the switching separatrix to be monotone, the proof of which is found in the Appendix.
\begin{thm}\label{thm:nec-suff}
Let the system~\eqref{eq:sys} satisfy Assumptions A1--A4. Then the following properties are equivalent:
\begin{enumerate}
\item If $\phi_f(t, s^0_f, \overline{\mu} h(\cdot,\overline{\tau}))$ belongs to $\cA(s^0_f)$ for all $t\ge 0$, then $\phi_f(t, s^0_f, \mu h(\cdot,\tau))$  belongs to $\cA(s^0_f)$ for all $t\ge 0$, and for all $\mu$, $\tau$ such that $0 <\mu\le \overline{\mu}$, $0 <\tau\le \overline{\tau}$.  
\item The set $\cS^-_f$ is simply connected. There exists a curve $\mu_f(\tau)$, which is a set of maximal elements of $\cS^-_f$ in the standard partial order. Moreover, the curve $\mu_f(\tau)$ is such that for any $\mu_1\in\mu_f(\tau_1)$ and $\mu_2\in\mu_f(\tau_2)$, $\mu_1\ge \mu_2$ for $\tau_1 < \tau_2$.
\end{enumerate}
\end{thm}

Theorem~\ref{thm:nec-suff} shows that the computation of the set $\cS^-_f$ is reduced to the computation of a curve $\mu_f(\tau)$. This result also provides a connection between the geometry of domains of attraction of the unforced system and the switching separatrix.  As shown next, Theorem~\ref{thm:nec-suff} can also be used to establish non-monotonicity of the switching separatrix.

\begin{rem}[Lorenz system revisited] Consider the Lorenz system from the previous section and three different pulses $u_{i}(t)=\mu_{i} h(t,\tau)$ with $\mu_{1}=24$, $\mu_{2}=25$, $\mu_{3}=26$, and $\tau=1$. Numerical solutions with increased accuracy show that the flows $\phi(t, s^0, u_1)$ and $\phi(t, s^0, u_3)$ converge to $s^0$, whereas $\phi(t, s^0, u_2)$ converges to $s^1$. Application of Theorem~\ref{thm:nec-suff} proves that the switching separatrix is not monotone.
\end{rem}

The major bottleneck in the direct application of Theorem~\ref{thm:nec-suff} is the verification of condition (1), which is generally computationally intractable. For example, condition (1) is satisfied if the partial order is preserved for control signals. That is for any $u\preceq_u v$, it should follow that $\phi_g(t, s^0_g, u)\preceq_x \phi_g(t, s^0_g, v)$ for all $t>0$. Although this property is weaker than monotonicity, it is not clear how to verify it. Monotonicity, on the other hand, is easy to check and implies condition (1) in Theorem~\ref{thm:nec-suff}. This is used in the following result.
\begin{thm}\label{thm:mon-switch-sep}
Let the system~\eqref{sys:low} satisfy Assumptions~A1--A4 and be monotone on $\cD_g\times\cU_\infty$. 
\begin{enumerate}
\item The set $\cS^-_g$ is simply connected. There exists a curve $\mu_g(\tau)$, which is a set of maximal elements of $\cS^-_g$ in the standard partial order. Moreover, the curve $\mu_g(\tau)$ is such that for any $\mu_1\in\mu_g(\tau_1)$ and $\mu_2\in\mu_g(\tau_2)$, $\mu_1\ge \mu_2$ for $\tau_1 < \tau_2$.
\item The set $\cS^+_g$ is simply connected. There exists a curve $\nu_g(\tau)$, which is a set of minimal elements of $\cS^+_g$ in the standard partial order. Moreover, the curve $\nu_g(\tau)$ is such that for any $\nu_1\in\nu_g(\tau_1)$ and $\nu_2\in\nu_g(\tau_2)$, $\nu_1\ge \nu_2$ for $\tau_1 < \tau_2$.
\item Let the system~\eqref{sys:low} be strongly monotone and $\partial \cA$ be the separatrix between the domains of attractions $\cA(s^0_f)$ and $\cA(s^1_f)$ of the unforced system~\eqref{sys:low}. Let additionally $\partial \cA$ be an unordered manifold, that is, there are no $x$, $y$ in $\partial\cA$ such that $x\succ_x y$. Then $\nu_g(\tau) = \mu_g(\tau)$ for all $\tau>0$ and the curve $\mu_g(\cdot)=\nu_g(\cdot)$ is a graph of a monotonically decreasing function.
\end{enumerate}
\end{thm}

We state implicitly in Theorem~\ref{thm:mon-switch-sep}, that if $\mu_g(\cdot) \ne \nu_g(\cdot)$, then the flow $\phi_g(t, s^0_g, \mu h(\cdot, \tau))$ does not converge to $s^0_g$ or $s^1_g$, since it may end up on the separatrix $\partial \cA$. We note that our computational procedure presented in Section~\ref{s:comp} does not require that $\mu_g(\tau) = \nu_g(\tau)$ or that $\mu_g(\cdot)$, $\nu_g(\cdot)$ are graphs of functions. Hence we treat point (3) in Theorem~\ref{thm:mon-switch-sep} as a strictly theoretical result, but remark that sufficient conditions for the separatrix $\partial \cA$ to be unordered are provided in Theorem~2.1 in~\cite{jiang2004saddle}. The most relevant condition to our case is that the flow of the unforced system is strongly monotone, which we also assume in Theorem~\ref{thm:mon-switch-sep}.

Besides $\mu_g(\cdot) \ne \nu_g(\cdot)$, there are other pathological cases. For example, applying constant input control signals $u = \mu$ typically results in a system~\eqref{eq:sys} with a different set of steady states than $s^0_f$ or $s^1_f$. Moreover, the number of equilibria may be different. Hence, with $\tau\rightarrow\infty$ the set $\cS^+_f$ typically does not contain the limiting control signal $u = \mu$. If the set of pairs $(\mu,\tau)$ resulting in these pathological cases is not measure zero, then the sets $\cl(\cS^+_f)$ and $\cl(\Rnn^2\backslash\cS^-_f)$ are not equal, which can complicate the computation of the switching sets. However, in many practical applications, the sets $\cl(\cS^+_f)$ and $\cl(\Rnn^2\backslash\cS^-_f)$ appear to be equal. Therefore in order to simplify the presentation we study only the properties of the set $\cS^-_f$. 

\begin{figure}[t]
\includegraphics[height=0.34\columnwidth]{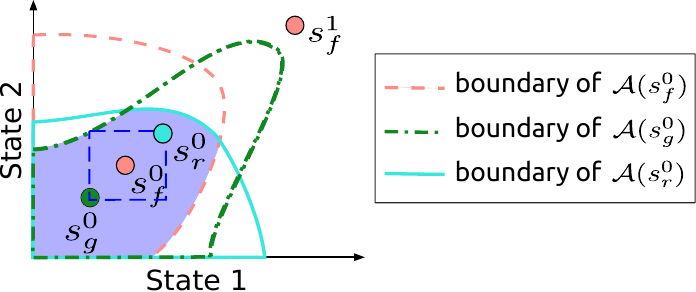}
\centering
\caption{A schematic depiction of the conditions~\eqref{eq:condss} and~\eqref{eq:condss2}. The condition~\eqref{eq:condss} ensures that all the steady states lie in the intersection of the corresponding domains of attractions (violet area). The steady state $s^1_f$ cannot lie in the dashed blue box due to condition~\eqref{eq:condss2}.}
\label{fig:propss}
\end{figure}

If the system $\dot x = f(x,u)$ to be controlled is not monotone, then the curve $\mu_f(\tau)$ may not be monotone, which is essential for our computational procedure. Instead, we estimate inner and outer bounds on the switching set provided that the vector field of the system can be bounded from above and below by vector fields of monotone systems. This is formally stated in the next result, while the proof is in the Appendix. 
\begin{thm}
\label{thm:comp-sys}
Let systems~\eqref{eq:sys},~\eqref{sys:low},~\eqref{sys:up} satisfy Assumptions~A1--A4. Let $\cD_M = \cD_g\cup\cD_f\cup\cD_r$, the systems~\eqref{sys:low} and~\eqref{sys:up} be monotone on $\cD_M \times \cU_\infty$ and 
\begin{align}
	g(x,u) \preceq_x f(x, u)\preceq_x r(x,u)\textrm{ on } \cD_M\times\cU.
\end{align}
Additionally assume that the stable steady states $s^0_g$,  $s^0_f$, $s^0_r$, $s^1_f$ satisfy 
\begin{gather}\label{eq:condss} 
s^0_g, s^0_f, s^0_r \in\inter\left(\cA(s^0_g)\cap\cA(s^0_f)\cap \cA(s^0_r)\right), \\
\label{eq:condss2} s^1_f\not \in \left\{z| s^0_g \preceq_x z \preceq_x s^0_r \right\}.
\end{gather}
Then the following relations hold:
\begin{equation}
\label{st:bound} \cS^-_{g} \supseteq \cS^-_{f} \supseteq \cS^-_{r}. 
\end{equation} 
\end{thm}
The technical conditions in~\eqref{eq:condss},~\eqref{eq:condss2} are crucial to the proof and are generally easy to satisfy. An illustration of these conditions is provided in Figure~\ref{fig:propss}.  Checking the condition \eqref{eq:condss2} reduces to the computation of the stable steady states, as does checking the condition~\eqref{eq:condss}. Indeed, to verify that $s^0_f$ belongs to the intersection of $\cA(s^0_g)$, $\cA(s^0_f)$, $\cA(s^0_r)$, we check if the trajectories of the systems~\eqref{sys:low},~\eqref{sys:up} initialised at $s^0_f$ with $u=0$ converge to $s^0_g$ and $s^0_r$, respectively, which is done by numerical integration. The computation of stable steady states can be done using the methods from~\cite{zwolak2004finding}.

In some applications, we need to find a subset of the pairs $(\mu,\tau)$ that switch the system~\eqref{eq:sys} from $s^0_f$ to $s^1_f$. Due to the inclusion $\cS^-_{g} \supseteq \cS^-_{f}$, existence of the system~\eqref{sys:low} allows to do that. In this case, we are only interested in finding the system~\eqref{sys:low}, hence the condition~\eqref{eq:condss2} is not required and the condition~\eqref{eq:condss} is transformed to $s^0_g$,  $s^0_f\in\inter\left(\cA(s^0_g)\cap\cA(s^0_f)\right)$.

\begin{rem}
The proofs of Theorems~\ref{thm:mon-switch-sep} and~\ref{thm:comp-sys} are adapted in a straightforward manner to the case when systems are monotone with respect to an order $\succeq_x$ induced by an arbitrary cone $K_x$, and the order $\succeq_u$ induced by $\Rnn$. In examples, however, we always assume that $K_x$ is an orthant, since it is hard to check monotonicity of the system with respect to an arbitrary cone.
\end{rem}

Theorem~\ref{thm:comp-sys} also provides a way of estimating the switching set under parametric uncertainty, which is stated in the next corollary. 
\begin{cor}
\label{cor:curve-unc}
Consider a family of systems $\dot x = f(x, u, p)$ with a vector of parameters $p$ taking values from a compact set $\cP$. Let the systems $\dot x = f(x, u, p)$ satisfy Assumptions~A1--A4 for every $p$ in $\cP$. Assume there exist parameter values $a$, $b$ in $\cP$ such that the systems $\dot x = f(x,u, a)$ and $\dot x = f(x,u, b)$ are monotone on $\cD_M\times \cU_\infty$, where $\cD_M=\mathop{\cup}\limits_{q\in\cP}\cD_{f(\cdot, \cdot, q)}$ and
\begin{align}
f(x,u, a)\preceq_x f(x, u, p) \preceq_x f(x, u, b),
\end{align}
for all $(x,u,p)\in\cD_M  \times \cU\times \cP$. Let also 
\begin{gather}
\label{eq:condssp} s_{f(\cdot, \cdot, p)}^0\in\inter\left(\mathop{\cap}\limits_{q\in\cP}\cA(s_{f(\cdot, \cdot, q)}^0)\right),\\
\label{eq:condssp2} s_{f(\cdot, \cdot, p)}^1 \not \in \left\{z| s^0_{f(\cdot, \cdot, a)} \preceq_x z \preceq_x s^0_{f(\cdot, \cdot, b)} \right\},
\end{gather}
for all $p$ in $\cP$. Then the following relation holds: 
\begin{equation}\label{eq:boundsparam}
 \cS^-_{f(\cdot,\cdot, a)} \supseteq \cS^-_{f(\cdot,\cdot, p)} \supseteq\cS^-_{f(\cdot,\cdot, b)} \quad\forall p\in\cP. 
\end{equation} 
\end{cor}
The proof follows by setting $g(x,u)= f(x,u, a)$ and $r(x,u)= f(x,u,b)$ and noting that the conditions in~\eqref{eq:condssp},~\eqref{eq:condssp2} imply the conditions in~\eqref{eq:condss},~\eqref{eq:condss2} in the premise of Theorem~\ref{thm:comp-sys}.

Theorem~\ref{thm:comp-sys} states that if the bounding systems~(\ref{sys:low}), (\ref{sys:up}) can be found, the switching sets $\cS^-_{g}$, $\cS^-_{r}$ can be estimated, thereby providing approximations on the switching set $\cS^-_{f}$. In what follows we provide a procedure to find monotone bounding systems if the system~\eqref{eq:sys} is near-monotone, meaning that by removing some interactions between the states the system becomes monotone (see~\cite{sontag2007monotone} for the discussion on near-monotone systems).

Let there exist a single interaction which is not compatible with monotonicity with respect to an order induced by $\Rnn^n$, and let this interaction be between the states $x_i$ and $x_j$. This happens if, for example, the $(i,j)$-th entry  in the Jacobian $\left\{\frac{\partial f_i}{\partial x_j}\right\}_{i,j}$ is smaller or equal to zero. A monotone system can be obtained by replacing the variable $x_j$ with a constant in the function $f_i(x,u)$, which removes the interaction between the states $x_i$ and $x_j$. If the set $\cD$ is bounded then clearly we can find $\overline x_j$ and $\underline x_j$ such that $\overline x_j \ge x_j \ge \underline x_j$ for all $x\in\cD$. If the set $\cD$ is not bounded, then we need to estimate the bounds on the intersection of $\cA(s^0_f)$ and the reachability set starting at $s^0_f$ for all admissible pulses. Let $g_k=r_k=f_k$ for all $k\ne i$, $g_i(x,u) = f_i(x,u)\bigl|_{x_j = \underline x_j}$, and $r_i(x,u) = f_i(x,u)\bigl|_{x_j = \overline x_j}$. It is straightforward to show that $\dot x = g(x,u)$, and $\dot x = r(x,u)$ are monotone systems and their vector fields are bounding the vector field $f$ from below and above, respectively. Note that in order to apply Theorem~\ref{thm:comp-sys} we still need to check if these bounding systems satisfy Assumptions A1--A4. 

In the case of Corollary~\ref{cor:curve-unc}, the procedure is quite similar. If the system $\dot x = f(x,u,p)$ is monotone for all parameter values $p$, then we can find $a$, $b$ if there exists a partial order in the parameter space. That is a relation $\preceq_p$ such that for parameter values $p_1$ and $p_2$ satisfying $p_1 \preceq_p p_2$ we have that
\[
f(x,u,p_1)\preceq_x f(x,u,p_2)~\forall x\in \cD,u\in\cU.
\]
If a partial order is found, the values $a$ and $b$ are computed as minimal and maximal elements of $\cP$ in the partial order $\preceq_p$. This idea is equivalent to treating parameters $p$ as inputs and showing that the system $\dot x = f(x,u,p)$ is monotone with respect to inputs $u$ and $p$.
\section{Computation of the Switching Separatrix\label{s:comp}}
The theoretical results in Section~\ref{s:theory} guarantee the existence of the switching separatrix for monotone systems, but in order to compute $\mu(\tau)$ we resort to numerical algorithms.

Given a pair $(\mu,\tau)$ we can check if this pair is switching the system using simulations (that is, numerically integrating the corresponding differential equation). If the curve $\mu(\tau)$ is a monotone function, then for every $\tau$ there exists a unique pulse magnitude $\mu=\mu(\tau)$. 
Let $\cT = \{\tau_i\}_{i=1}^N$ be such that $\tau_{\rm min} = \tau_1\le\tau_i \le \tau_{i+1}\le\tau_N=\tau_{\rm max}$ for all $i$. Clearly, for every $\tau_i$ we can compute the corresponding $\mu_i$ using bisection. We start the algorithm by computing the value $\mu_1$ corresponding to $\tau_1$. Due to monotonicity of the switching separatrix, the minimal switching magnitude $\mu_2$ for the pulse length $\tau_2$ is smaller or equal to $\mu_1$. Therefore, we can save some computational effort by setting the upper bound on the computation of $\mu_2$ equal to $\mu_1$. The computation of the pairs $(\mu,\tau)$ can be parallelised by setting the same upper bound on $\mu_{i}$, $\cdots$, $\mu_{i + N_{\rm par}}$, where $N_{\rm par}$ is the number of independent computations. The procedure is summarised in Algorithm~\ref{alg:u-min}, where $\cM_{\rm min}$ and $\cM_{\rm max}$ are the sets of pairs $(\mu,\tau)$ approximating the switching separatrix from below and above, respectively. 

\begin{algorithm}[t]
\caption{Bisection Algorithm for Computation of the Switching Separatrix}
\label{alg:u-min}
\begin{algorithmic}[1]
\State {\bf Inputs:} The system $\dot x= f(x, u)$ with initial state $s_f^0$, final state $s_f^1$, tolerance $\varepsilon$, simulation time $t_e$, a grid $\cT=\{\tau_i\}_{i = 1}^N$, an upper bound on the magnitude $\mu^{\rm up}$, the number of used processors $N_{\rm par}$.
\State {\bf Outputs:} finite sets $\cM^{\rm min}$ and $\cM^{\rm max}$
     \For{$i= 1,\dots, [N/N_{\rm par}]$}
     \State Set $\mu^l_j = 0$, $\mu_j^{u} = \mu^{\rm up}$, 
     \For{$j = 1,\dots, N_{\rm par}$}
     \While{$\mu^u_j - \mu_j^l>\varepsilon$}
            \State $\mu_j^c = (\mu_j^u+\mu_j^l)/2$
            \State Set $u^c = \mu_j^c h(\cdot,\tau_{i(N_{\rm par} -1)+j})$
            \If{$\phi_f(t_e; s_f^0, u^c) \preceq_x s_f^1$} 
              \State $\mu_j^c = \mu_j^l$
            \Else
              \State $\mu_j^c = \mu_j^u$
            \EndIf
     \EndWhile
     \EndFor
     \State $\cM^{\rm min} = [\cM^{\rm min}, \mu_1^l,\dots \mu^l_{N_{\rm par}}]$
     \State $\cM^{\rm max} = [\cM^{\rm max}, \mu_1^u,\dots \mu^u_{N_{\rm par}}]$
     \EndFor
\end{algorithmic}
\end{algorithm}

\begin{figure}[b]
  \centering
  \includegraphics[width = 0.7\columnwidth]{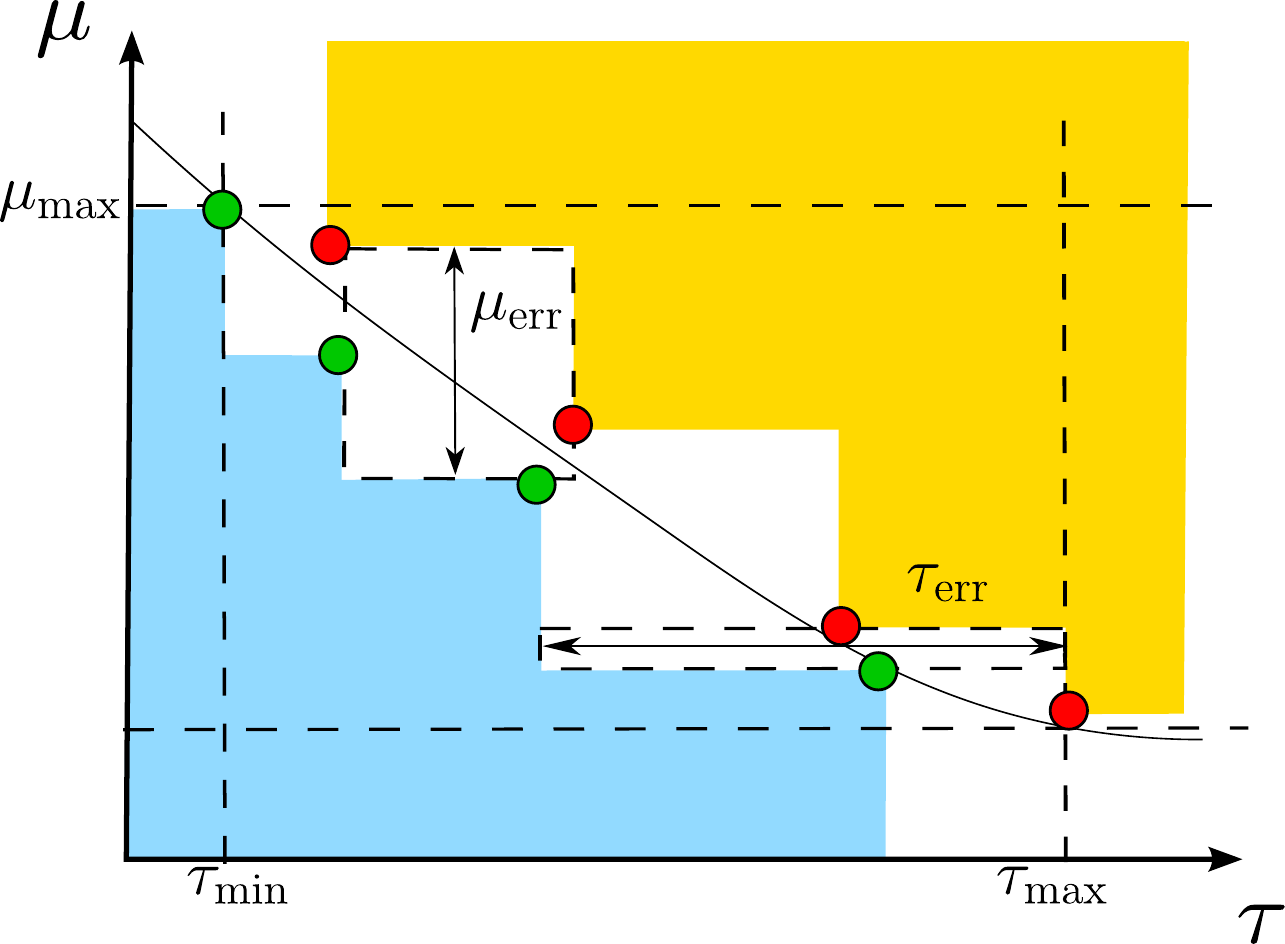}
  \caption{Illustration of the error of computation of the switching separatrix between the values $\tau_{\rm min}$, $\tau_{\rm max}$.  The black curve is the switching separatrix to be computed, the red and green circles are the upper and lower bounding points, respectively. The switching separatrix should lie between the coloured regions due to its monotonicity. The values $\mu_{\rm err}$ and $\tau_{\rm err}$ are the largest height and width of boxes inscribed between the coloured regions, respectively.} 
  \label{fig:err-plot}
\end{figure}
In order to evaluate the error of computing the switching separatrix consider Figure~\ref{fig:err-plot}. According to the definitions in the caption of Figure~\ref{fig:err-plot} we define the relative error of the approximation as
\[E_{\rm rel} = (\mu_{\rm err}/(\mu_{\rm max}-\mu_{\rm min})+\tau_{\rm err}/(\tau_{\rm max}-\tau_{\rm min}))/2.\]
 Note that, even if the green and red circles lie very close to each other the relative error can be substantial. In numerical simulations we use a logarithmic grid for $\tau$, which yields a significantly lower relative error in comparison with an equidistant grid. This can be explained by an observation that in many numerical examples $\mu(\tau)$ appears to be an exponentially decreasing curve.

\begin{algorithm}[t]
\caption{Computation of Switching Separatrix Based on Random Sampling}
\label{alg:rand-comp}
\begin{algorithmic}[1]
\State {\bf Inputs:} The system $\dot x= f(x, u)$ with initial state $s_f^0$, final state $s_f^1$, total number of samples $N$, simulation time $t_e$, lower and upper bounds on $\tau$, $\tau_{\rm min}$ and  $\tau_{\rm max}$ respectively, the numbers $N_{\rm gr}$, $N_{\varepsilon}$, probability distribution $\delta$
\State {\bf Outputs:} sets $\cM^{\rm min}$ and $\cM^{\rm max}$
     \State Compute $\mu_{\rm min}$ and $\mu_{\rm max}$ using bisection for values $\tau_{\rm min}$ and $\tau_{\rm max}$
     \State Set $N_{\rm par} = 2 (N_{\rm gr}+N_{\varepsilon})$
     \State Set $\cM^{\rm max}=\cM^{\rm min} = \{(\mu_{\rm max},\tau_{\rm min}), (\mu_{\rm min},\tau_{\rm max})\}$
     \For{$i= 1,\dots, [N/N_{\rm par}]$}
	 \State Compute the values $\mu_{\rm err}$, $\tau_{\rm err}$, and the corresponding boxes $\cB_\mu$, $\cB_\tau$.
	 \State Generate $N_{\rm gr}$ samples $(\mu,\tau)$ in each of the boxes $\cB_\mu$ and $\cB_\tau$ using a probability distribution $\delta$
	 \State Generate randomly $2 N_{\varepsilon}$ samples
     \For{$j = 1,\dots, N_{\rm par}$}
			\State Check if the samples $(\mu,\tau)$ switch the system
     \EndFor
     \State Update and prune the sets $\cM^{\rm min}$, $\cM^{\rm max}$ 
     \EndFor
\end{algorithmic}
\end{algorithm}
There are a few drawbacks in Algorithm~\ref{alg:u-min}. Firstly, it requires a large number of samples. Secondly, the choice of the grid is not automatic, which implies that for switching separatrices with different geometry the relative error on the same grid may be drastically different. Finally, the algorithm relies on the assumption that $\mu(\tau)$ is a graph of a monotone function, which may not be true. In order to overcome these difficulties, we have derived Algorithm~\ref{alg:rand-comp} based on random sampling, which converges faster than Algorithm~\ref{alg:u-min}, has higher sample efficiency, does not require a predefined grid and the graph assumption. Some of the steps in Algorithm~\ref{alg:rand-comp} require additional explanation:

Step 7. Find two boxes: the box $\cB_\mu$ with the maximal height (denoted as $\mu_{\rm err}$) and the box $\cB_\tau$ with the maximal width (denoted as $\tau_{\rm err}$) that can be inscribed between the coloured regions as depicted in Figure~\ref{fig:err-plot}. 

Step 9. Generate $N_{\varepsilon}$ samples of $\tau$ using a probability distribution $\delta$ between $\tau_{\rm min}$ and $\tau_{\rm max}$. For every $\tau$ generate a value $\mu$ using a distribution $\delta$ such that $\mu$ lies in the area between the coloured regions. Repeat this step by first generating $\mu$ between $\mu_{\rm min}$ and $\mu_{\rm max}$ using a distribution $\delta$, and then generating $\tau$ for every generated $\mu$ in the area between the coloured regions. 

Step 13. First, we update the sets $\cM^{\rm min}$, $\cM^{\rm max}$ by adding the samples that do not switch and switch the system, respectively. 
Then if there exist two pairs $(\mu_1,\tau_1)$ and $(\mu_2,\tau_2)$ in the set $\cM^{\rm min}$ (resp., $\cM^{\rm max}$) such that $\mu_1 \le \mu_2$ and $\tau_1\le\tau_2$, then delete the pair $(\mu_1, \tau_1)$ from the set $\cM^{\rm min}$ (resp., the pair $(\mu_2,\tau_2)$ from the set $\cM^{\rm max}$). 

Note that Step 11 is the most computationally expensive part of the algorithm and its computation is distributed into $N_{\rm par}$ independent tasks. In our implementation, we chose $\delta$ as a Beta distribution with parameters $1$ and $3$ and adjusted the support to a specific interval. Note that the set between the coloured regions is getting smaller with every generated sample, hence the relative error of Algorithm~\ref{alg:rand-comp} is a non-increasing function of the total number of samples. In fact, numerical experiments show that this function is on average exponentially decreasing. After the sets $\cM^{\rm min}$ and $\cM^{\rm max}$ are generated one can employ machine learning algorithms to build a closed form approximation of a switching separatrix (e.g., Sparse Bayesian Learning~\cite{tipping2001sparse}; see also \cite{wipf2008new}, \cite{wei2014admm} for efficient algorithms). 

\begin{figure}[t]
  \centering
  \includegraphics[width = 0.88\columnwidth]{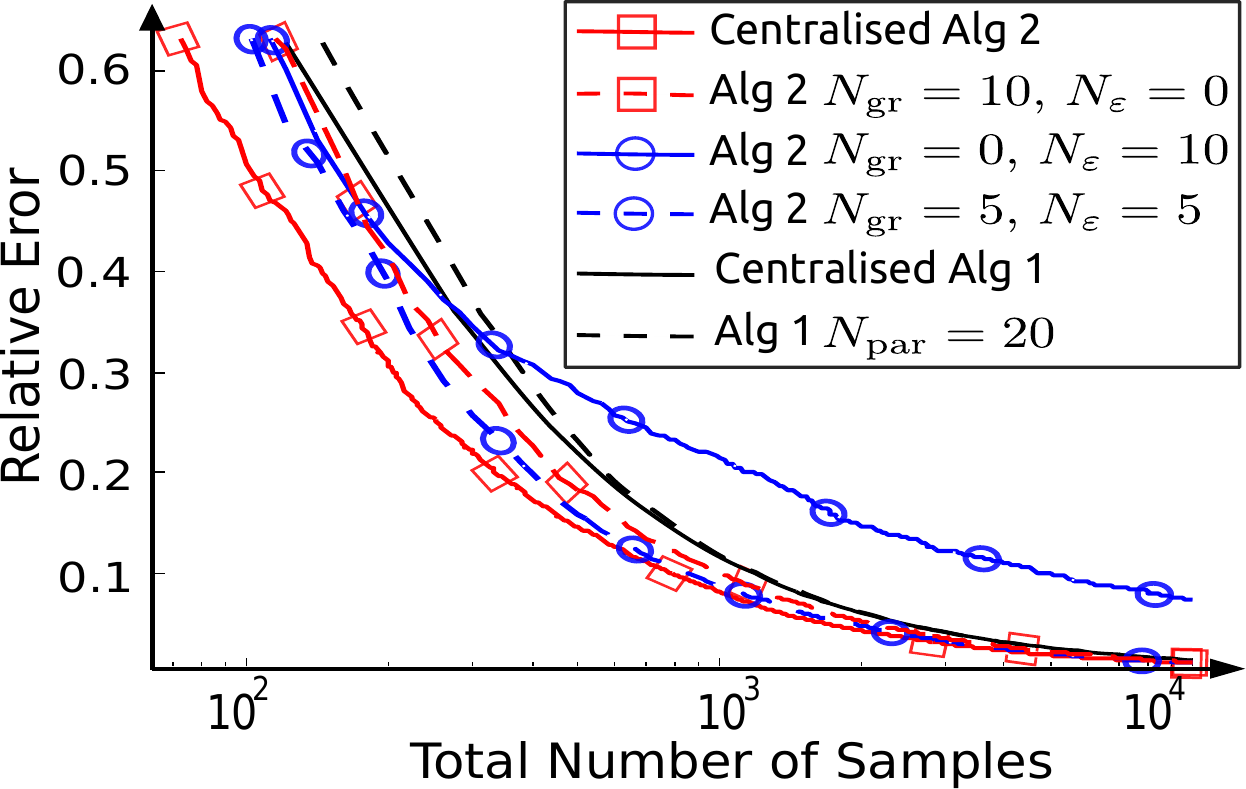}
  \caption{Average error against total number of generated samples. The curves corresponding to   Algorithm~\ref{alg:u-min} are computed by a single run of the algorithm. The curves corresponding to $N_{\varepsilon} = 0$ are averages over ten runs of Algorithm~\ref{alg:rand-comp}, while the curves for $N_{\varepsilon} > 0$ are the averages over twenty runs of Algorithm~\ref{alg:rand-comp}. Recall that $N_{\rm par} = 2(N_{\rm gr} + N_{\varepsilon})$ for Algorithm~\ref{alg:rand-comp}.}
  \label{fig:err}
\end{figure}
\section{Illustration of Theoretical Results on the LacI-TetR Switch\label{s:laci-tetr}}

\subsection{Evaluation of the Computational Algorithm\label{ss:eff}}
Here we will compare Algorithm~\ref{alg:u-min} and~\ref{alg:rand-comp} with different parameter values, as well as their distributed implementations on the LacI-TetR switch introduced in Subsection~\ref{ss:mech}. Note that Algorithm~\ref{alg:rand-comp} does not depend explicitly on the dynamics of the underlying system, but depends only on the generated pairs $(\mu,\tau)$. Therefore, the convergence and sample efficiency results presented here will be valid for a broad class of systems. In Figure~\ref{fig:err}, we compare the error against the total number of generated samples. Since checking if a sample switches the system is the most expensive part of both algorithms, the total number of samples reflects the computational complexity.  In the case of Algorithm~\ref{alg:rand-comp} with $N_{\varepsilon} = 0$ the randomisation level is not high, hence an average over ten runs is sufficient to demonstrate the average behaviour of this algorithm. Note that both curves corresponding to Algorithm~\ref{alg:rand-comp} with $N_{\varepsilon}=0$ outperform the curves corresponding to Algorithm~\ref{alg:u-min} in terms of accuracy. 

\begin{table}[t]\centering
 \caption{Sample efficiency $N_{\rm eff}$ in percent. In the notation $x\pm y$, $x$ stands for the emperical mean, and $y$ for the emperical standard deviation. Recall that $N_{\rm par} = 2(N_{\rm gr} + N_{\varepsilon})$ for Algorithm~\ref{alg:rand-comp}.\label{tab:se}}
  \medskip
 \begin{tabular} {l c }
 Algorithm    & $N_{\rm eff}$ \\
\hline
\hline
 Alg.~\ref{alg:u-min}  with $N_{\rm par} = 1$                           & $14\%$  \\
 Alg.~\ref{alg:u-min}  with $N_{\rm par} = 20$                          & $13.3\%$   \\
 Alg.~\ref{alg:rand-comp}  with $N_{\rm gr} = 1$, $N_{\varepsilon} = 0$  & $28.21\pm 0.72\%$\\
 Alg.~\ref{alg:rand-comp}  with $N_{\rm gr} = 10$, $N_{\varepsilon} = 0$ & $19.90\pm 0.45\%$\\
 Alg.~\ref{alg:rand-comp}  with $N_{\rm gr} = 5$, $N_{\varepsilon} = 5$  & $39.68\pm 0.76\%$\\
  Alg.~\ref{alg:rand-comp}  with $N_{\rm gr} = 0$, $N_{\varepsilon} = 10$& $49.36\pm 0.55\%$\\  
 \hline
\end{tabular} 
\end{table}

Some computational effort in Algorithm~\ref{alg:rand-comp} goes into computing the error. However, this effort appears to be negligible in comparison with numerically solving a differential equation for a given pair $(\mu,\tau)$ even for such a small system as the toggle switch. We run the simulations on a computer equipped with Intel Core i7-4500U processor and 8GB of RAM. Using the centralised version of Algorithm~\ref{alg:rand-comp} we achieved on average a relative error equal to $0.0448$ in $87.65$ seconds, while it took $89.17$ seconds to obtain a relative error equal to $0.0842$ with Algorithm~\ref{alg:u-min}. Note that for systems with a larger number of states the difference may be larger. 

In Table~\ref{tab:se}, we compare the sample efficiency of Algorithms~\ref{alg:u-min} and~\ref{alg:rand-comp} with different input parameters, which we define as
\[
N_{\rm eff} = |\cM^{\rm min}\cup\cM^{\rm max}|/N
\]
where $N$ is the total number of generated samples, and $|\cM^{\rm min}\cup\cM^{\rm max}|$ is the number of samples in the set $\cM^{\rm min}\cup\cM^{\rm max}$. Results in Table~\ref{tab:se} indicate that Algorithm~\ref{alg:rand-comp} has higher sample efficiency than Algorithm~\ref{alg:u-min}. 

Our results also indicate that Algorithm~\ref{alg:rand-comp} with $N_{\rm gr} = 5$, $N_{\varepsilon} = 5$ has on average a higher empirical convergence rate and a higher sample efficiency than Algorithm~\ref{alg:rand-comp} with $N_{\rm gr} = 10$, $N_{\varepsilon} = 0$. This indicates that combination of non-zero $N_{\rm gr}$, $N_{\varepsilon}$ improves convergence and sample efficiency, which can be explained as follows. When the total number of generated samples is low, we do not have sufficient information on the behaviour of the switching separatrix. Therefore we need to explore this behaviour by randomly generating samples, before we start minimising the relative error. This idea is similar to the so-called exploration/exploitation trade-off in reinforcement learning~\cite{busoniu2010reinforcement}. 
\begin{table}[t]\centering
 \caption{Parameter values for systems in Subsection~\ref{ss:ss_ts}. The unspecified parameter values are the same as in~\eqref{eq:ts-par-val}.}\label{tab:systems}
  \medskip
  \begin{tabular} {l c c}
  							& $\cF^i_{\rm upper}$      & $\cF^i_{\rm lower}$ \\
\hline
\hline
\multirow{2}{*}{$i = 1$}    & $p_1 = 40$, $p_4 = 0.05$ & $p_1 = 20$, $p_4 = 0.01$, \\
                            & $p_6 = 30$, $p_9 = 0.1$ &$p_6 = 45$, $p_9 = 0.3$,\\
\hline
\multirow{3}{*}{$i = 2$}   & $p_1 = 40$, $p_4 = 0.05$ & $p_1 = 20$, $p_4 = 0.01$, \\
                            & $p_6 = 30$, $p_9 = 0.1$ &$p_6 = 45$, $p_9 = 0.3$,\\ 
                            & $p_2 =4$, $p_7 =1$      & $p_2 =1$, $p_7= 4$, \\
\hline
\multirow{3}{*}{$i = 3$}   & $p_1 = 40$, $p_4 = 0.05$ & $p_1 = 20$, $p_4 = 0.01$, \\
                           & $p_6 = 30$, $p_9 = 0.1$ &$p_6 = 45$, $p_9 = 0.3$,\\
                           & $p_5 =1$, $p_{10} = 2$  & $p_5 =3$, $p_{10} = 1$,\\  
\hline
\end{tabular} 
\end{table}
\subsection{Switching in the LacI-TetR System \label{ss:ss_ts}}
\begin{figure}[t]
  \centering
  \includegraphics[width = 0.7\columnwidth]{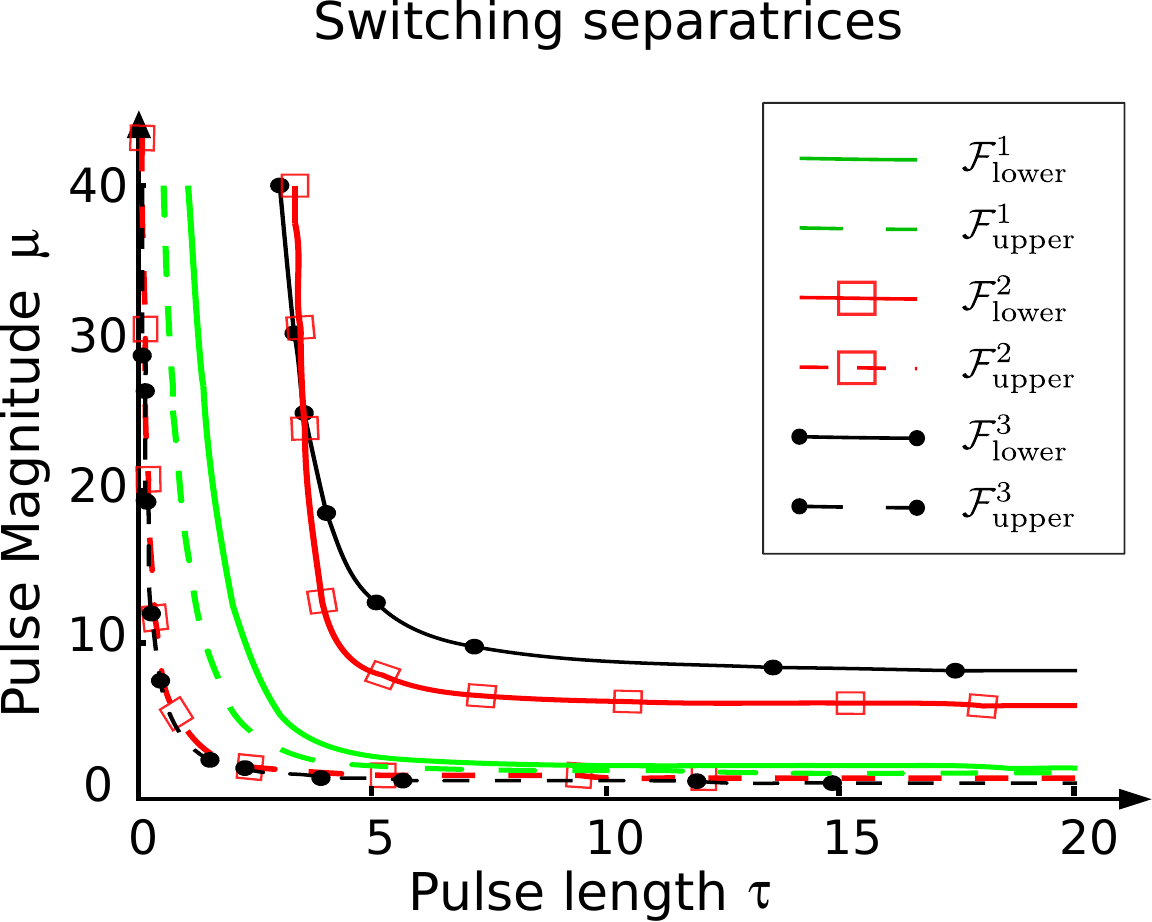}
  \caption{Switching separatrices for the LacI-TetR system~\eqref{eq:ts-par}. The system parameters are defined in Table~\ref{tab:systems}.}
  \label{fig:ts-separ}
\end{figure}

\emph{Robust Switching in the LacI-TetR System.} In Table~\ref{tab:systems}, we specify different versions of the LacI-TetR system by varying some of the parameters in~\eqref{eq:ts-par-val}. After that we compute the switching separatrices and plot them in Figure~\ref{fig:ts-separ}. Note that the separatrices for $\cF^2_{\rm lower}$, $\cF^3_{\rm lower}$ intersect, which happens since the vector fields $f_2$ and $f_3$ of systems $\cF^2_{\rm lower}$ and $\cF^3_{\rm lower}$, correspondingly, are not comparable. This means that there exists a set $\cX = \{(x,u) \in \cD \times \cU\}$ on which $f_2(x,u) \not \preceq_x f_3(x,u)$ and $f_3(x,u) \not \preceq_x f_2(x,u)$. 

The separatrices for $\cF^1_{\rm upper}$, $\cF^1_{\rm lower}$ lie very close to each other despite the number of parameters varied and the level of variations. This is not true for the separatrices $\cF^2_{\rm upper}$, $\cF^2_{\rm lower}$ and $\cF^3_{\rm upper}$, $\cF^3_{\rm lower}$. This indicates that the switching separatrix is sensitive to variations of some parameters more than to variations of other. In our case, this happens because the variations in parameters $p_5$, $p_{10}$, $p_2$, and $p_7$ affect significantly the positions of the stable steady states. Therefore, pulses with significantly smaller magnitudes are required to switch the systems $\cF^2_{\rm lower}$ and $\cF^3_{\rm lower}$ in comparison with $\cF^2_{\rm upper}$ and $\cF^3_{\rm upper}$, respectively. 

\emph{Switching in LacI-TetR Systems with Perturbed Dynamics.} Consider the following three-state system
  \begin{equation}
    \label{eq:non-mon}
    \cF = \left\{\begin{array}{rl}
      \dot x_1 &= \dfrac{1000}{1 + x_3^{2}} - 0.4 x_1, \\
      \dot x_2 &= \dfrac{1000}{1 + x_1^{4}}  - 4 x_2 + u, \\
      \dot x_3 &= p_{1}  + p_{2}x_1 + p_{3} \frac{x_1}{x_1+1}+ 5 x_2  - 0.3 x_3. 
    \end{array}\right. 
  \end{equation}
Consider two nominal systems $\cF^1$ and $\cF^2$ specified in Table~\ref{tab:nm-systems} by changing parameter values for $p_1$, $p_{2}$, $p_{3}$. In Table~\ref{tab:nm-systems}, the notations $\cG_{\rm upper}^i$ and  $\cG_{\rm lower}^i$ stand for the upper and lower bounding system of $\cF^i$ (for $i=1,2$).
\begin{table}[b]\centering
 \caption{Parameter values for systems with different parameter variations, where $z_0^i(1)$ is the first component of the initial point $z_0^i$ of the system $\cF^i$ for $i$ equal to $1$ and $2$. \label{tab:nm-systems}}
  \medskip
  \begin{tabular} {c c c}
\hline
 $\cF^1$ & $\cG_{\rm lower}^1$ &$\cG_{\rm upper}^1$ \\
 $p_{2} = 0.1$    & $p_1= 0$    & $p_1 = 0.1 z_0^1(1)$\\
 $p_1= p_{3}=0$   &  $p_{2} = p_{3}=0$ & $p_{2} =p_{3}=0$ \\
 \hline
 $\cF^2$ & $\cG_{\rm lower}^2$ &$\cG_{\rm upper}^2$ \\
 $p_{3}=0.1$      &  $p_1= 0$ & $p_1 =0.1 \frac{z_0^2(1)}{z_0^2(1)+1}$\\
 $p_1= p_{2} = 0$ & $p_{2} = p_{3}=0$ &  $p_{2} = p_{3}=0$ \\
\hline
\end{tabular} 
\end{table}

\begin{figure}[t]
  \centering
\includegraphics[width = 0.7\columnwidth]{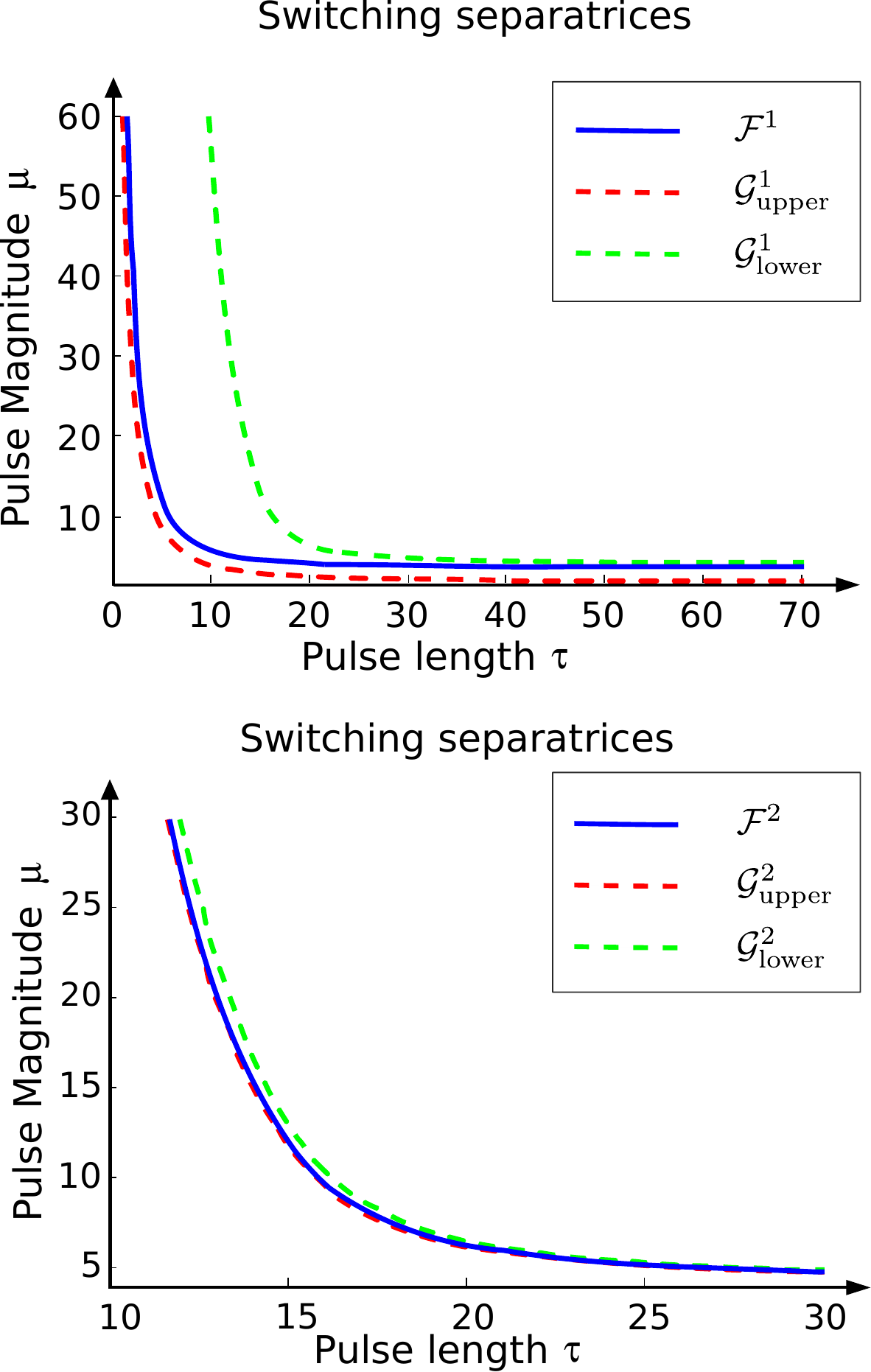}
  \caption{The switching separatrices for systems $\cF^1$ (the upper panel) and $\cF^2$ (the lower panel) with lower and upper bounds. }
  \label{fig:nm-separ}
\end{figure}

Consider first the system $\cF^1$ and compute the bounding systems $\cG_{\rm upper}^1$ and  $\cG_{\rm lower}^1$. We will follow the procedure described in Section~\ref{s:theory}. Using Proposition~\ref{prop:kamke}, it is easy to check that with a positive value for $p_{2}$ this system is not monotone with respect to any orthant. Hence, we need to bound the term $p_{2}x_1$ by constants in order to obtain monotone bounding systems. By simulating the system we observe that $x_1$ lies in a bounded interval between $0$ and $z_0^1(1)$, where $z_0^1(1)$ is the first component of the initial point $z_0^1$. Hence, we can build an upper $\cG_{\rm upper}^1$ and a lower $\cG_{\rm lower}^1$ bounding system for the nominal one $\cF^1$. We take the system $\cG_{\rm upper}^1$  with the same parameter values as the nominal one except for $p_{2}$, which is equal to zero, and $p_1$ equal to $0.1 z_0^1(1)$. Similarly, we choose the system $\cG_{\rm lower}^1$ with $p_{2} = 0$, and $p_1 = 0$.  The results can be seen in the upper panel of Figure~\ref{fig:nm-separ}. We have computed the approximation of the switching separatrix for $\cF^1$ by simulating the system on a mesh grid in $\mu$ and $\tau$. This requires a considerable computational effort, but provides a decent approximation. Note that the switching separatix for $\cF^1$ appears to be a monotone curve, even though it cannot be guaranteed. However, this can be guaranteed for the separatrices of the bounding systems, which are monotone.

Now let us compute the bounds on the switching separatrix of the nominal system $\cF^2$, where the Michaelis-Menten term ($x_1/(x_1 + 1)$) prevents the system from being monotone. In a similar fashion as for the case of $\cF^1$, we can build an upper $\cG_{\rm upper}^2$ and a lower $\cG_{\rm lower}^2$ bounding systems for $\cF^2$. This results in the switching separatrices depicted in the lower panel of Figure~\ref{fig:nm-separ}. The bounds on the switching separatrix for $\cF^2$ are tighter in comparison with the bounds of the switching separatrix for $\cF^1$. Note that in the case of the system $\cF^1$, we use the following bound $0 \le x_1 \le z^1_0(1)$, while in the case of the system $\cF^2$, we use the bound $0 \le \frac{x_1}{x_1+1} \le \frac{z_0^2(1)}{z_0^2(1)+1}$, where $z_0^2(1)$ is the first component of the initial point of the system $\cF^2$. At the same time the numbers $z_0^1(1)$, $z_0^2(1)$ are of order $O(10^3)$, hence clearly $\frac{z_0^2(1)}{z_0^2(1)+1}\ll z_0^1(1)$. This means that the bounds on the vector field of $\cF^2$ are tighter than the bounds on the vector field of $\cF^1$. This in turn indicates that tighter bounds on the vector field entail tighter bounds on the separatrix. 

\section{Further Counterexamples and Applications\label{s:examples}}
\begin{figure}[t]\centering
  \includegraphics[width = 0.7\columnwidth]{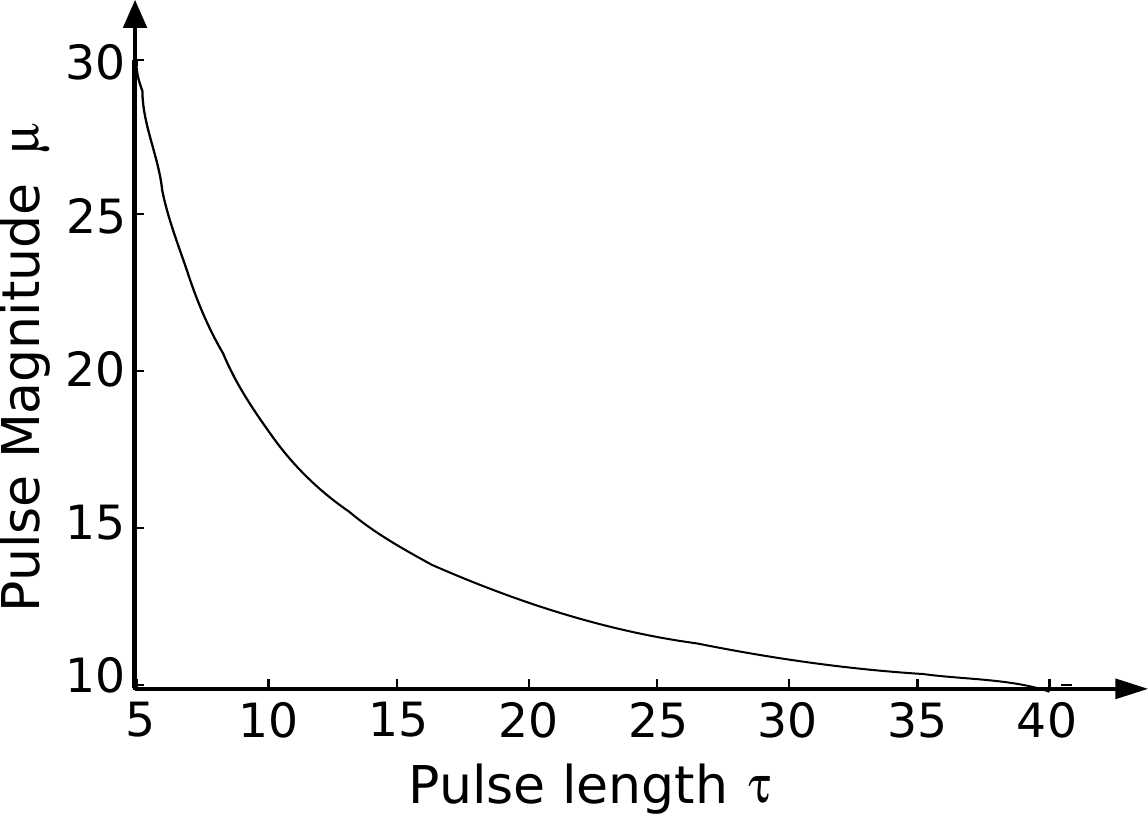}
\caption{Approximation of the switching separatrix for the toxin-antitoxin system. } \label{fig:ta}
\end{figure}

\emph{Toxin-Antitoxin System.} Consider a model of a toxin-antitoxin system studied in~\cite{cataudella2013conditional}
\begin{align*}
\dot T &= \frac{\sigma_T}{(1 + \frac{[A_f][T_f]}{K_0})(1+\beta_M [T_f])} - \frac{1}{(1+\beta_C [T_f])} T \\
\dot A &= \frac{\sigma_A}{(1 + \frac{[A_f][T_f]}{K_0})(1+\beta_M [T_f])} - \Gamma_A  A + u\\
\varepsilon \dot{[A_f]} &= A - \left([A_f] + \frac{[A_f] [T_f]}{K_T} + \frac{[A_f] [T_f]^2}{K_T K_{T T}}\right) \\ 
\varepsilon \dot{[T_f]} &= T - \left([T_f] + \frac{[A_f] [T_f]}{K_T} + 2\frac{[A_f] [T_f]^2}{K_T K_{T T}}\right), 
\end{align*}
where $A$ and $T$ is the total number of toxin and antitoxin proteins, respectively, while $[A_f]$, $[T_f]$ is the number of free toxin and antitoxin proteins. In~\cite{cataudella2013conditional}, the authors considered the model with $\varepsilon = 0$. In order to simplify our analysis we set $\varepsilon = 10^{-6}$. If the parameters are chosen as follows:
\begin{gather*}
\sigma_T = 166.28,~~K_0 = 1,~~\beta_M = \beta_c =0.16,~~\sigma_A = 10^2\\
\Gamma_A = 0.2,~~K_T = K_{TT} = 0.3;
\end{gather*}
the system is bistable with two stable nodes. But the system is not monotone and we were not able to find bounding systems satisfying Assumptions A1-A4. Nevertheless, we approximated the switching sets and the switching separatrix on a mesh grid in $\mu\in[10, 30]$ and $\tau\in[5,40]$. We had $200$ points in the grid for $\tau$ and $80$ points in the grid for $\mu$. In Figure~\ref{fig:ta} we plot a curve which separates the pairs $(\mu,\tau)$ which toggle the system from those which do not. As the reader may notice the curve appears to be monotone. We can provide some intuition behind this phenomenon. With $\varepsilon$ tending to zero, we can apply singular perturbation theory (cf.~\cite{KhalilNonlinearControl}) to eliminate the states $[A_f]$, $[T_f]$. Numerical computations indicate that the reduced order system is not monotone in $\Rnn^2$, however, it is monotone around the stable equilibria, which may explain monotonicity of the switching separatrix of the full order system.

\emph{Switching in a Mass Action Kinetics System.} The following system was considered in~\cite{wilhelm2009smallest} 
\begin{align*}
\dot x_1 &= f_1(x_1, x_2)= 2 k_1 x_2 - k_2 x_1^2 - k_3 x_1 x_2 - k_4 x_1 + \beta u, \\
\dot x_2 &= f_2(x_1, x_2)= k_2 x_1^2 - k_1 x_2.
\end{align*}
We can assume without loss of generality that $k_2 = 1$, since we can remove one of the parameters using a simple change of variables. In~\cite{wilhelm2009smallest} it is shown that the unforced system has three equilibria:
\begin{gather*}
s^0 = \begin{pmatrix}0 \\ 0 \end{pmatrix}~
s^u = \begin{pmatrix}\frac{k_1 - \sqrt{k_1 L}}{2 k_3}\\ \left(\frac{\sqrt{k_1} - \sqrt{L}}{2 k_3}\right)^2\end{pmatrix}~
s^1 = \begin{pmatrix}\frac{k_1 + \sqrt{k_1 L}}{2 k_3}\\ \left(\frac{\sqrt{k_1} + \sqrt{L}}{2 k_3 }\right)^2 \end{pmatrix}
\end{gather*}
where $L = k_1 - 4 k_3 k_4$, $s^0$, $s^1$ are stable nodes, and $s^u$ is a saddle. Naturally the system is bistable if $L>0$, globally asymptotically stable if $L < 0$ and a saddle-node bifurcation occurs if $L = 0$. The system is not monotone in the positive orthant, since the derivative of $f_1(x_1, x_2)$ with respect to $x_2$ is equal to $2 k_1  - k_3 x_1$ and it is negative for $x_1> 2 k_1/k_3$. However, we are interested in pairs $(\mu,\tau)$, which are close to the separatrix and affect the system in the domain of attraction of $s^0$. It can be verified that all equilibria $s^0$, $s^u$ and $s^1$ lie inside the box $\cD = \{ x_1, x_2 | 0 \le x_1 \le 2 k_1/k_3\}$, hence the system is monotone inside $\cD$ and non-monotonicity outside of this set does affect our problem. 

The derivatives of $f_1$, $f_2$ with respect to $k_1$ and $k_2$ do not have the same sign hence the system is not monotone with respect to parameters $k_1$ and $k_2$. These terms appear due to so-called mass action kinetics, which are common in biological applications and hence this problem is met often. We have already removed the parameter $k_2$ from the consideration, which simplifies the problem. A straightforward solution is to treat every instance of $k_1$ as an independent parameter. Hence we  have a vector of parameters $[k_{11},~k_3,~k_4,~k_{12}]$, where $k_{11}$ is the instance of $k_1$ entering the first equation, and  $k_{12}$ is the instance of $k_1$ entering the second equation. Let $k_1 \in [7.7,~8.3]$, $k_3 \in [1,~1.2]$, $k_4 \in [1,~1.2]$ and consider the lower bounding parameter vector $p_l = [7.7,~1.2,~1.2,~8.3]$, and the upper bounding parameter vector $p_u = [8.3,~1,~1,~7.7]$. We apply Corollary~\ref{cor:curve-unc} only to relatively small perturbations in parameters, since with larger variations the system becomes mono- or unstable. There is no indication that this problem is unique to this system, and does not appear in other mass-action systems.

\begin{figure}[t]\centering
  \includegraphics[width = 0.7\columnwidth]{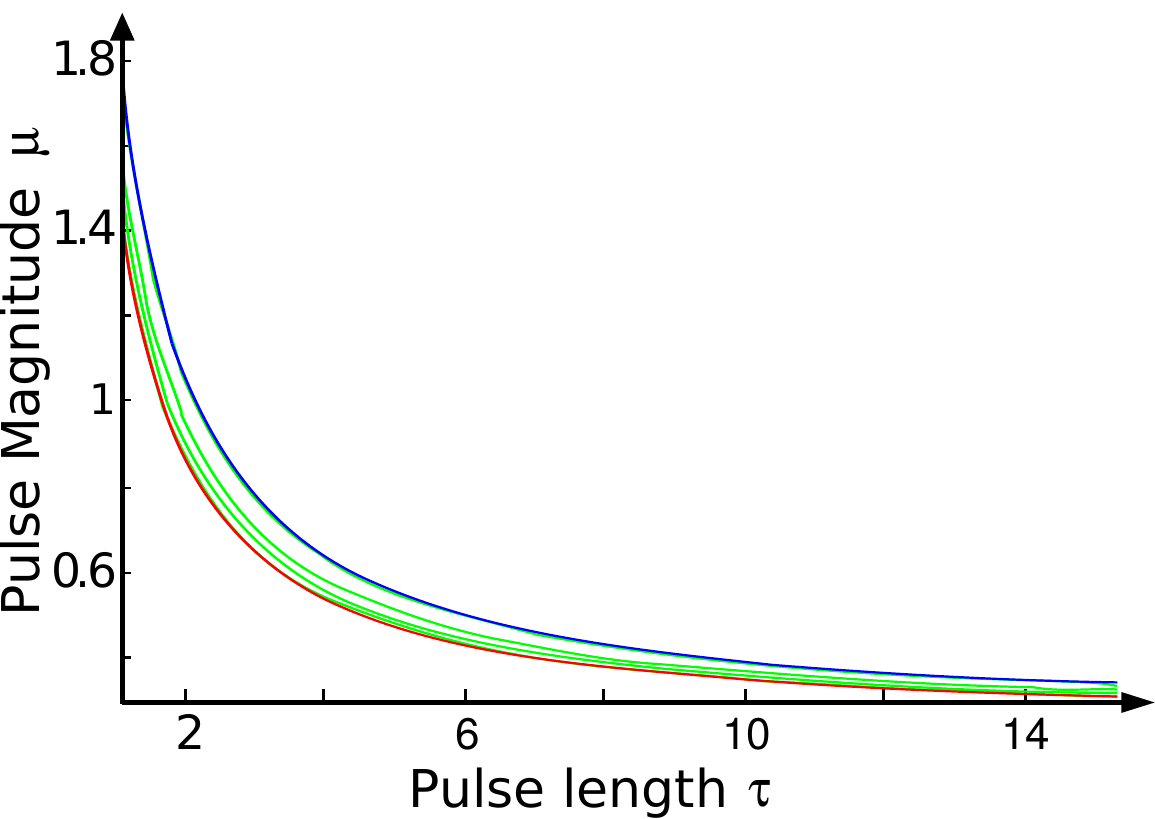}
\caption{Switching separatrices for a mass action kinetics system with different values for $k_1$. The blue curve is the switching separatrix for $k_1 = 6$, the red curve is the switching separatrix for $k_1 = 10$. The green curves are the switching separatrices for $k_1 = 6.1$, $7.5$, $8.5$, $9.9$.} \label{fig:param-sweep}
\end{figure}

We conclude this example by performing a sweep for the parameter $k_1\in [6,10]$, while $k_2 = k_3 = 1$ (see Figure~\ref{fig:param-sweep}). Numerical simulations suggest that for any $k_1 \in (6, 10)$ the switching separatrix appears to lie between the blue and the red curves, which are switching separatrices for $k_1 = 6$ and $k_1 = 10$, respectively. Again we can only provide some intuition behind this observation. It is straightforward to verify that the gradient of $s^u$ with respect to $k_1$ is a negative vector, and the gradient of $s^1$ with respect to $k_1$ is a positive vector. Hence the equilibria depend on $k_1$ in the way which is consistent with a behaviour of a monotone system. This example indicates that the behaviour of the equilibria may be one of the necessary conditions allowing the switching separatrix to be a monotone curve and change monotonically with respect to parameter variations.

\begin{figure*}
  \centering
 \includegraphics[height = 0.45\columnwidth]{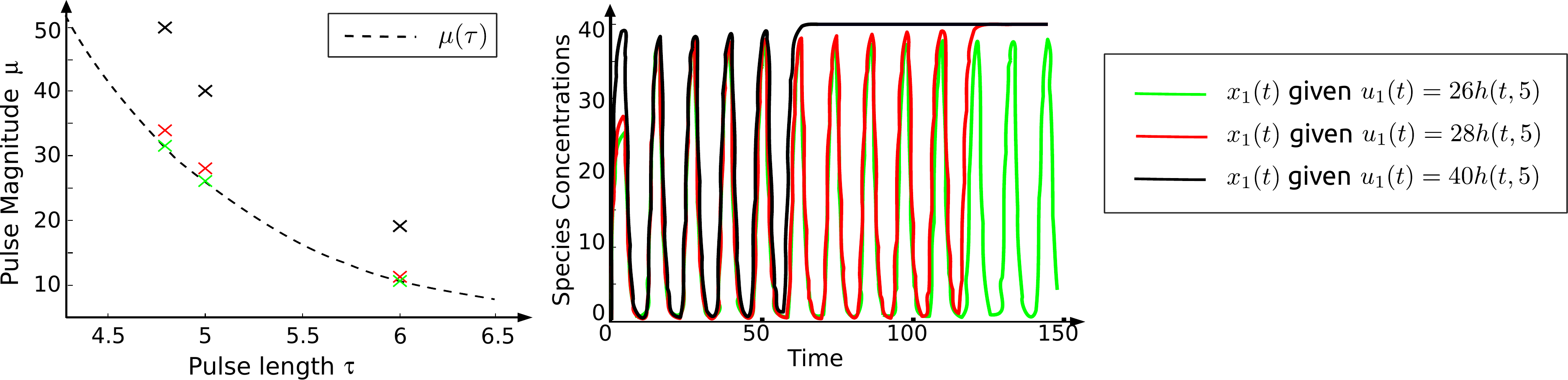}
  \caption{Switching between steady states in a generalised repressilator system. All trajectories generated by the pairs $(\mu,\tau)$ corresponding to the black crosses in the left panel will converge to a steady state with the same rate as the black curve in the right panel. Similar correspondence is valid for the red and green crosses in the left panel and the red and green curves in the right panel. This observation indicates that the closer the pair $(\mu,\tau)$ is to the switching separatrix the longer oscillations will persist. }
  \label{fig:gr8-sw}
  \end{figure*}

\emph{Shaping Pulses to Induce Oscillations in an Eight Species Generalised Repressilator.} An eight species generalised repressilator is an academic example, where each of the species represses another species in a ring topology. The corresponding dynamic equations for a symmetric generalised repressilator are as follows: 
  \begin{equation}
    \label{eq:gr}
    \begin{aligned}
   \dot x_1 &= \frac{p_{1}}{1 + (x_{8}/p_2)^{p_3}} + p_4 - p_5 x_1 + u_1, \\
   \dot x_2 &= \frac{p_{1}}{1 + (x_{1}/p_2)^{p_3}} + p_4 - p_5 x_2 + u_2, \\
   \dot x_i &= \frac{p_{1}}{1 + (x_{i-1}/p_2)^{p_3}} + p_4 - p_5 x_i,~\forall i = 3,\dots 8,  
    \end{aligned}
  \end{equation}
where $p_1 = 40$, $p_2 = 1$, $p_3 = 3$, $p_4 = 0.5$, and $p_5 = 1$. This system has two stable equilibria $s^1$ and $s^2$ and is monotone with the respect to $P_x \R^8\times P_u\R^2$, where $P_x= \diag([1,~-1,~1,~-1,~1,~-1,~1,~-1])$, $P_u = \diag([1,~-1])$. It can actually be shown that the system is strongly monotone in the interior of $\Rnn^8$ for all positive parameter values. The control signal $u_1$ can switch the system from the state $s^1$ to the state $s^2$, while the control signal $u_2$ can switch the system from the state $s^2$ to the state $s^1$. The switching separatrix for the control signal $u_1$ is depicted in the left panel of Figure~\ref{fig:gr8-sw}. Note that the separatrix is identical for the control signal $u_2$, since the repressilator is symmetric.

Numerical simulations suggest that the trajectories exhibit an oscillatory behaviour, while switching between the stable steady states using a pulse. This is in agreement with previous studies that showed the existence of unstable periodic orbits~\cite{Strelkowa10} in a generalised repressilator. Switching trajectories of species $x_1$ for various pairs $(\mu,\tau)$ are depicted in the right panel of Figure~\ref{fig:gr8-sw}. The observations made in the caption of Figure~\ref{fig:gr8-sw} indicate that the closer the pair $(\mu,\tau)$ is to the switching separatrix the longer oscillations will persist.

We can set up another control problem: \emph{to induce oscillations in the generalised repressilator}. One can address the problem by forcing the trajectories to be close to the unstable periodic orbit of the system, which, however, is very hard to compute. In~\cite{sootla2013tracking}, it was proposed to track other periodic trajectories instead. However, the solution was very computationally expensive and offering little insight into the problem. Here we will use pulses to induce oscillations as was proposed in~\cite{Strelkowa10}. In contrast to~\cite{Strelkowa10}, we provide a way to shape all possible pulses inducing oscillations. 

Let the initial point be $s^1$. We can shape the control signal $u_1$ to switch to the state $s^2$. Once we have reached an $\varepsilon$-ball around the state $s^2$, we can shape the control signal $u_2$ to switch back to the state $s^1$ and so on. During switching we will observe oscillations depending on the position of the pair $(\mu,\tau)$ with respect to the switching separatrix. Now we need to define an automatic way of switching between the steady states. Let $\cM$ be equal to $\{z \bigl| s^1 \preceq_x z \preceq_x s^2\}$. It can be verified that the unstable equilibrium lies in $\cM$, which typically holds for monotone systems. Moreover, the trajectories observed in Figure~\ref{fig:gr8-sw} lie in $\cM$ due to monotonicity. Let $\varepsilon > 0$ and $\cM_\varepsilon=\{z \bigl| s^1 + \varepsilon P_x \bfone \preceq_x z \preceq_x s^2- \varepsilon P_x \bfone\}$, where $\bfone$ is the vector of ones. Clearly $\cM_\varepsilon\subset \cM$ and if $\varepsilon$ is small enough then oscillating trajectories lie in $\cM_\varepsilon$. Since the repressilator is symmetric we can assume that the shape of pulses for both $u_1$ and $u_2$ is the same. In this case we can formalise our control strategy as follows. If the event $x(t_e) \preceq_x s^1 + \varepsilon P_x \bfone$ occurs at time $t_e$, then
\begin{align*}
&u_1(\cdot) = \mu h(\cdot,t_e + \tau) \quad
&u_2(\cdot) = 0  
\end{align*}
\begin{figure}[b]
\centering
  \includegraphics[height = 0.45\columnwidth]{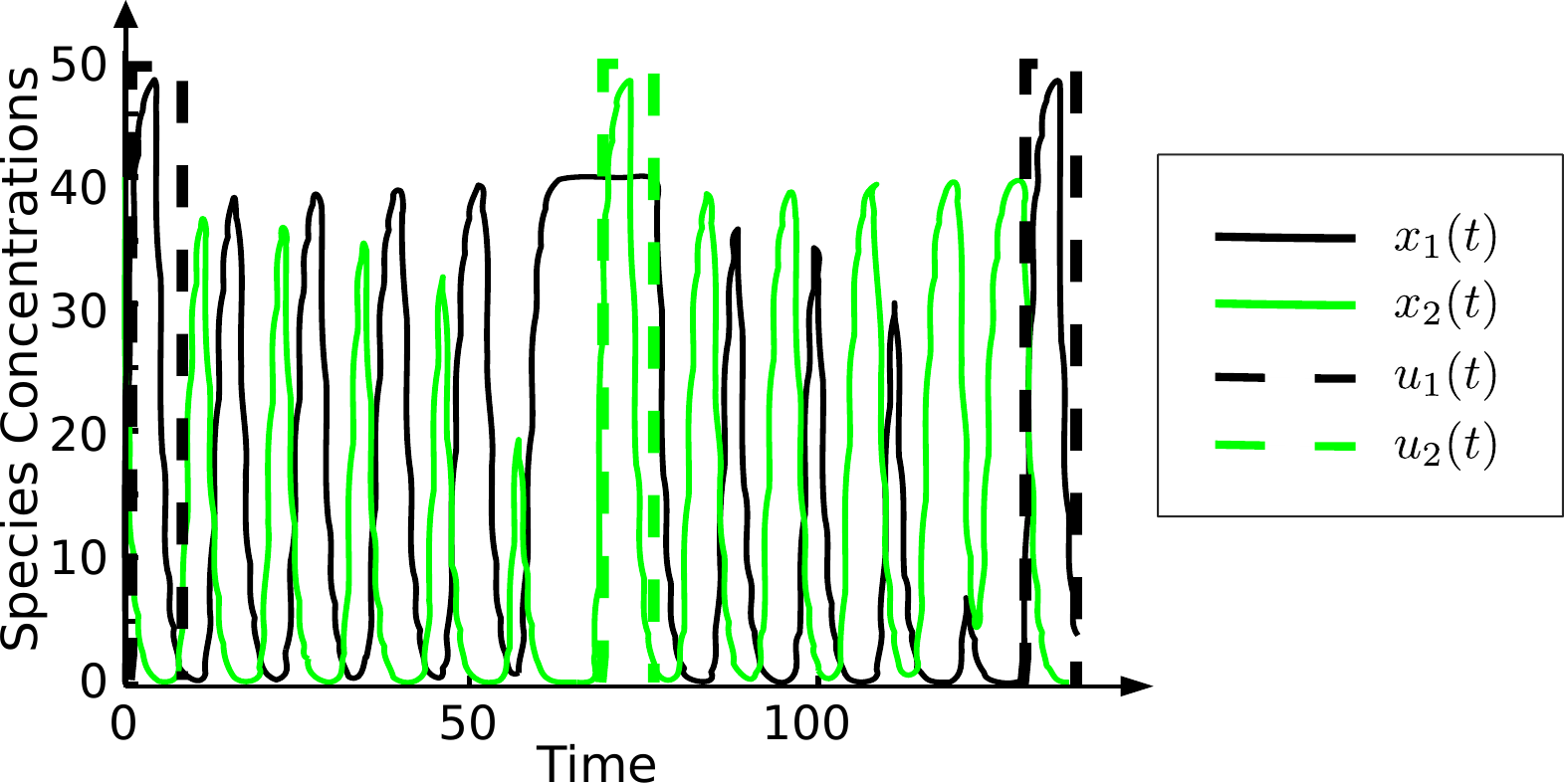}
\caption{Inducing oscillatory behaviour in the generalised repressilator system with eight species. The pulses for both $u_1$ and $u_2$ are equal, and are generated using a pair $(\mu,\tau) = (48, 4.8)$. The pair $(48, 4.8)$ lies relatively far from the switching separatrix, hence the time between switches is not large.} \label{fig:gr8-osc}
\end{figure}
If the event $x(t_e)\succeq_x s^2 - \varepsilon P_x \bfone$ occurs at time $t_e$, then 
\begin{align*}
&u_1(\cdot) = 0 \quad
&u_2(\cdot) = \mu h(\cdot,t_e +\tau)
\end{align*}
Note that we change the entire control signals when the event occurs, which we assume happens at some time $t_e$. Due to this fact, the pulse $\mu h(\cdot,t_e + \tau)$ is of length $\tau$. The resulting trajectories for the species $x_1$ and $x_2$, as well as control signals are depicted in Figure~\ref{fig:gr8-osc}. Our control algorithm falls into the class of event-based control, with the events occurring if $x(t_e)$ leaves $\cM_\varepsilon$. For any small enough $\varepsilon$, our control strategy will induce oscillations. 

\section{Conclusion and Discussion}

In this paper we have presented a framework for shaping pulses to control bistable systems. Our main motivation comes from control problems arising in Synthetic Biology, but the results hold in other classes of bistable systems.
We considered the problem of switching between stable steady states using temporal pulses. We showed that the problem is feasible, if the flow of the controlled system can be bounded from above and below by 
flows of monotone systems. We presented a detailed analysis of the conditions needed for switching, together with an algorithm to compute the pulse's length and duration. We illustrated the theory with a number of case studies and counterexamples that shed light on the limitations of the approach and highlight the need for further theoretical tools to control bistable non-monotone systems.

Throughout this work we did not take into account stochasticity in the model dynamics, which can be particularly important in biochemical systems~\cite{elowitz2002stochastic}. Noisy bistable dynamics can be controlled, for example, using reinforcement learning algorithms as the ones described in~\cite{sootla2013tracking,sootla2014toggle}. These approaches, however, require large amounts of measurement data that are typically impractical to acquire. A promising extension to our results is the switching problem in stochastic bimodal systems. This requires the use of the so-called stochastically monotone Markov decision processes, for which a whole new set of theoretical tools needs to be developed. Work in this direction started in~\cite{Sootla2015PO} and the references within, addressing the extension of the concept of monotonicity to stochastic systems. 

\bibliography{Biblio}
\appendices
\section{Proofs}
\begin{proof-of}{Proposition~\ref{prop:bif}}
{\bf (1)} Here we simply need to notice that by monotonicity with $t\rightarrow\infty$ we have
\begin{gather*}
\xi(\mu)\leftarrow \phi_g(t, s^0_g, \mu) \preceq_x \phi_g(t, s^0_g, \lambda)\rightarrow \xi(\lambda).
\end{gather*}
 Similarly we can show that $\eta(\mu)\preceq_x \eta(\lambda)$. \\  
{\bf (2)} First, we need to show that $\xi(\mu)\in\cA(s^0_g)$ for all  $0<\mu<\mu_{\rm min}$. This is straightforward, since due the definition of $\mu_{\rm min}$ the flow $\phi_g(t, s^0_g, \mu h(\cdot, \tau))$ converges to $s^0_g$ for all the pairs $(\mu,\tau)\in \cS^-$ for  $0<\mu<\mu_{\rm min}$. Hence  the limit $\lim\limits_{t\rightarrow\infty}\phi_g(t, s^0_g, \mu)$ belongs to $\cA(s^0_g)$. \\  
Now, we show that $s^0_g\prec_x s^1_g$. Consider $u = 0$ and $v = \lambda h(\cdot,\tau)$  such that $(\lambda,\tau)\in \cS^+$. Therefore we have
\begin{gather*}
s^0_g = \phi_g(t, s^0_g, u) \preceq_x \phi_g(t, s^0_g, v) \rightarrow s^1_g,
\end{gather*}
with $t\rightarrow\infty$. Since $s^0_g$ is not equal to $s^1_g$, we have $s^0_g\prec_x s^1_g$. Now the claim $\xi(\mu) \prec_x \eta(\mu)$ for all $0<\mu<\mu_{\rm min}$ follows by monotonicity.\\
{\bf (3)} Consider $\mu = \mu_{\rm min} +\varepsilon$ and $\tau$ large enough that the pair $(\mu,\tau)\in\cS^+$. Hence the flow $\phi_g(t, s^0_g, \mu h(\cdot,\tau))$ converges to $s^1_g$. By monotonicity we have that $\phi_g(t, s^0_g, \mu)\succeq_x \phi_g(t, s^0_g, \mu h(\cdot,\tau))$, which implies that $\xi(\mu)\succeq_x s^1_g$ for arbitrarily small $\varepsilon>0$. Since $\xi(\mu_{\rm min} -\varepsilon)$ lies in $\cA(s^0_g)$ we have that  $\|\xi(\mu_{\rm min} -\varepsilon) - \xi(\mu)\|_2 > 0$, which proves the claim.
\end{proof-of}
\begin{proof-of}{Theorem~\ref{thm:nec-suff} }
{\bf $\mathbf{1)\Rightarrow 2)}$}~~{\bf A.} It is straightforward to verify that the premise of Theorem~\ref{thm:nec-suff} implies that any point lying in the set $\cS^-_f$ is path-wise connected to a point in the neighbourhood of the origin. In order to show that the set is simply connected, it is left to prove that there are no holes in the set $\cS^-_f$. Let $\eta(\mu,\tau)$ be a closed curve which lies in $\cS^-_f$. Consider the set 
\[
\cS^\eta = \left\{(\mu, \tau)\bigl| 0<\mu \le \mu^\eta,0<\tau\le\tau^\eta, (\mu^\eta,\tau^\eta)\in\eta(\mu,\tau)\right\}.
\]
Since the set $\cS^-_f$ is in $\Rp^2$, the set $\cS^\eta$ contains the set enclosed by the curve $\eta(\mu,\tau)$. It is straightforward to show that $\cS^\eta$ is a subset of $\cS^-_f$ by the premise of the theorem. Hence there are no holes in the area enclosed by the arbitrary curve $\eta\in\cS^-_f$. Since the curve $\eta$ is in $\R^2$ we can shrink this curve continuously to a point, which belongs to the set $\cS^-_f$. Since the curve is an arbitrary closed curve in $\cS^-_f$, the set $\cS^-_f$ is simply connected. \\
~~~{\bf B.} Let us show here that there exists a set of maximal elements in $\cS^-_f$. Let a pair $(\mu^u, \tau^u)$ not belong to $\cS^-_f$. If there exists a pair $(\mu,\tau)\in\cS^-_f$ such that $\mu \ge \mu^u$, $\tau \ge \tau^u$, then by the arguments above the pair $(\mu^u,\tau^u)$ must also belong to $\cS^-_f$. Hence, all pairs $(\mu, \tau)$ such that $\mu \ge \mu^u$, $\tau \ge \tau^u$ do not belong to $\cS^-_f$. This implies that there exists a set of maximal elements of $\cS^-_f$ in the standard partial order, which is a segment of the boundary of $\cS^-_f$ excluding the points with $\mu$ and $\tau$ equal to zero. \\
~~~{\bf C.} It is left to establish that the set of maximal elements is unordered. Let the mapping $\mu_f(\tau)$ denote the set of maximal elements of $\cS^-_f$ and let $\tau_1 < \tau_2$. Since the mapping $\mu_f(\tau)$ are the maximal elements in $\cS^-_f$, we cannot have $\mu_f(\tau_1) < \mu_f(\tau_2)$. Hence, $\mu_f(\tau_1) \ge \mu_f(\tau_2)$, for all $\tau_1 < \tau_2$. \\
{$\mathbf{2)\Rightarrow 1)}$} The claim follows directly from the fact that there exists a set of maximal elements $\mu_f(\tau)$ in the simply connected set $\cS_f^-$.
\end{proof-of}
\begin{proof-of}{Theorem~\ref{thm:mon-switch-sep}}
{\bf 1)} Due to Assumption A4, there exists at least one point $(\mu^l, \tau^l)$ in $\cS_g^-$. Let us show that if a pair $(\mu^l, \tau^l)$ belongs to $\cS^-_g$, then all pairs $(\mu,\tau)$ such that $0<\mu \le \mu^l$, $0<\tau \le \tau^l$ also belong to $\cS^-_g$.  By the definition of the order in $u$, for every $0<\mu \le  \mu^l$, $0<\tau \le \tau^l$  we have $0 \preceq_u\mu h(t,\tau)\preceq_u  \mu^l h(t, \tau^l)$. The following relation is then true
\[
s_g^0 \preceq_x \phi_{g}(t; s_g^0,\mu h(\cdot,\tau))\preceq_x \phi_{g}(t; s_g^0, \mu^l h(\cdot,\tau^l)).
\]
By assumption, there exists a $T$ such that for all $t>T$ the flow $\phi_{g}(t; s_g^0,\mu^l h(\cdot,\tau^l))$ belongs to $\cA(s_g^0)$ and converges to $s_g^0$. Therefore  $\phi_{g}(t; s_g^0,\mu h(\cdot,\tau))$ converges to $s_g^0$ with $t\rightarrow+\infty$, and consequently the pair $(\mu, \tau)$ does not toggle the system and thus belongs to $\cS^-_g$. Therefore, by Theorem~\ref{thm:nec-suff} $\mu_g(\tau_1) \ge \mu_g(\tau_2),$ for all $\tau_1 < \tau_2$. \\
{\bf 2)} Due to Assumption A4, there exists at least one point $(\mu^l, \tau^l)$ in $\cS_g^+$. Similarly to point {\bf 1)} above, we can show that, if a pair $(\mu^l, \tau^l)$ belongs to $\cS^+_g$, then by continuity of solutions to~\eqref{sys:low} there exist $\overline{\varepsilon}>0$, $\overline{\delta}>0$ such that the pairs $(\mu+\varepsilon,\tau+\delta)$ also belong to $\cS^+_g$ for all $0<\varepsilon<\overline{\varepsilon}$, $0<\delta<\overline{\delta}$. Hence the set $\cS^+_g$ has a non-empty interior. The rest of the proof is the same as the proof of the implication {\bf $1)\Rightarrow 2)$} in~Theorem~\ref{thm:nec-suff}. \\
{\bf 3)} We prove the result by contradiction. Let there exist a $\tau$ and an interval $I = (\mu_1, \mu_2)$ such that for all $\mu \in I$ the flow $\phi_g(t, s^0_g, \mu h(\cdot, \tau))$ does not converge to $s^0_g$ or $s^1_g$, but belongs to the interior of $\cD_g$. This means that the flow $\phi_g(t, s^0_g, \mu h(\cdot, \tau))$ evolves on the separatrix $\partial \cA$ between domains of attraction $\cA(s^0_g)$ and $\cA(s^1_g)$ for all $t>\tau$. Let $\mu_1$, $\mu_2$ belong to $I$ and $\mu_1 < \mu_2$, which implies that $\phi_g(t, s^0_g, \mu_1 h(\cdot, \tau))\ll_x \phi(t, s^0_g, \mu_2 h(\cdot, \tau))$ and both flows belong to $\partial \cA$. This in turn implies that the set $\partial \cA$ contains comparable points, that is, the set $\partial \cA$ is not unordered. We arrive at a contradiction, and hence the interval $I$ is empty and for any $\tau$ there exists a unique $\mu_g(\tau)$. This is equivalent to $\mu_g(\cdot)$ being a graph of a function. Using similar arguments, we can show that the inverse mapping $\mu_g^{-1}(\tau)$ is also a graph of a function, which indicates that $\mu_g(\cdot)$ is a decreasing function.\\
Similarly, we can show that for any $\mu$ the minimum value of $\tau_2 -\tau_1$, such that the pairs $(\mu,\tau_1 -\varepsilon)\in\cS_g^-$, $(\mu,\tau_2 +\varepsilon)\in\cS_g^+$ $\forall \varepsilon>0$, is equal zero. This readily implies that $\mu_g(\tau) = \nu_g(\tau)$ and completes the proof.
\end{proof-of}
Before we proceed with the proof of Theorem~\ref{thm:comp-sys} we will need two additional results: one is the so-called comparison principle for control systems and the other is concerned with geometric properties of the regions of attractions of monotone systems. 
\begin{figure}[t]
\includegraphics[width=0.4\columnwidth]{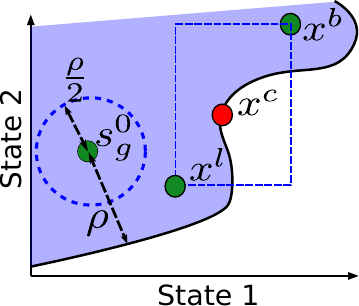}
\centering
\caption{An illustration to the proof of Lemma~\ref{lem:box-mon} for a two-state system. We assume that $x^b$, $x^l$ lie in $\cA(s^0_g)$ (violet area) and $x^b\succeq_x x^c \succeq_x x^l$ with $x^c$ lying on the boundary of $\partial\cA(s^0_g)$. We show that, if the trajectory $\phi_g(t,x^c,0)$ is on the boundary of $\cA(s^0_g)$, it has to converge to $s^0_g$, which cannot be true due to monotonicity of the system.}
\label{fig:propra}
\end{figure}
\begin{lem}\label{lem:box-mon}
Let the system $\dot{x} = g(x,0)$ satisfy Assumption~A1 and be monotone on $\cA(s^0_g)$, where $s^0_g$ is a stable steady state and $\cA(s^0_g)$ is its domain of attraction. Let $x^b$ and $x^l$ belong to $\cA(s^0_g)$. Then all points $z$ such that $x^l \preceq_x z \preceq_x x^{b}$ belong to $\cA(s^0_g)$.
\end{lem}
\begin{pf}
We will show the result by contradiction. Let $x^l$, $x^b$ belong to $\cA(s_g^{0})$, let $x^c$ be such that $x^l \preceq_x x^c \preceq_x x^b$ and not belong to $\cA(s^0_g)$. Without loss of generality assume that $x^c$ belongs to the boundary of $\cA(s^0_g)$ (see Figure~\ref{fig:propra}). Therefore the flow  $\phi_g(t, x^c, 0)$ is on the boundary of $\cA(s^0_g)$. Let the distance between $s_g^0$ and this boundary be equal to $\rho$. Clearly there exists a time $T_1$ such that for all $t>T_1$ the following inequalities hold
\begin{align*}
\|s^0_g - \phi_g(t, x^b, 0)\|_2 < \rho/2,\\
\|s^0_g - \phi_g(t, x^l, 0)\|_2 < \rho/2.
\end{align*}
Moreover, there exists a time $T_2>T_1$ such that for all $t>T_2$ and all $z$ such that 
$
\phi_g(t, x^l, 0) \preceq_x z \preceq_x \phi_g(t, x^b,0)
$ 
we have $\|s^0_g - z\|_2 < \rho/2$. Now build a sequence $\{x^n\}_{n=1}^{\infty}$
converging to $x^c$ such that all $x^n$ lie in $\cA(s^0_g)$ and $x^l\preceq_x x^n \preceq_x x^b$. Due to monotonicity on $\cA(s^0_g)$, we have 
\[
\phi_g(t, x^l, 0) \preceq_x \phi_g(t, x^n, 0) \preceq_x \phi_g(t, x^b, 0)
\]
for all $n$ and $t$. Hence, for all $t>T_2$, we also have that $\|s^0_g - \phi_g(t, x^n, 0)\|_2 < \rho/2$.
Since the sequence $\{x^n\}_{n=1}^{\infty}$ converges to $x^c$, by continuity of solutions to~\eqref{sys:low}, for all $t>T_2$ we have $\|s^0_g - \phi_g(t, x^c, 0)\|_2 \le \rho/2$, which is a contradiction since $\|s^0_g - \phi_g(t, x^c, 0)\|_2 \ge \rho$ for all $t$.
\end{pf}
\begin{lem}
\label{lem:comp-prin-control}
 Consider the dynamical systems $\dot x = f(x, u)$ and $\dot x = g(x, u)$ satisfying Assumption~A1. Let one of the systems be monotone on $\cD_M\times\cU_{\infty}$. If $g(x,u) \succeq_x f(x,u)$ for all $(x, u)\in\cD_M\times\cU$ then for all $t$, and for all $x_2\succeq_x x_1$, $u_2\succeq_u u_1$ we have $\phi_{g}(t;x_2, u_2) \succeq_x \phi_{f}(t; x_1, u_1)$.
\end{lem}
\begin{pf}
Without loss of generality let $\dot x = g(x, u)$ be monotone with respect to $\Rnn^n$. Let $\bfone$ be a vector of ones, $x_2^m = x_2 + 1/m \cdot\bfone$, and $\dot x = g(x,u)+\bfone/m$. Denote the flow of this system $\phi_m(t; x_2^m, u_2)$. Clearly for a sufficiently small $t$ the condition $\phi_m(t; x_2^m, u_2) \gg_x \phi_{f}(t; x_1, u_1)$ holds. Assume there exists a time $s$, for which this condition is violated. That means that for some $i$ we have $\phi_m^i(t; x_2^m, u_2) > \phi_{f}^i(t; x_1, u_1)$  for all $0\le t<s$, where the superscript $i$ denotes an $i$-th element of the vector. While at time $s$ we have $\phi_m^i(s; x_2^m, u_2) = \phi_{f}^i(s; x_1, u_1)$. Hence we conclude that  
\begin{gather}
\frac{d }{dt}(\phi_m^i(t; x_2^{m}, u_2)- \phi_{f}^i(t; x_1, u_1))\Bigl|_{t=s} \le 0. \label{cond:contr}
\end{gather}
However,
\begin{gather}
\frac{d \phi_{f}^i(t; x_1, u_1)}{dt}\Bigl|_{t=s}  = f_i(\phi_{f}(s; x_1, u_1), u_1) < \label{cond:contr-lem-comp-1}\\
g_i(\phi_{f}(s; x_1, u_1), u_1)+ 1/m  \le \label{cond:contr-lem-comp-2}\\
g_i(\phi_{m}(s; x_2^m, u_2), u_2)+ 1/m = \frac{d \phi_{m}^i(t; x_2^m, u_2)}{dt} \Bigl|_{t=s}.\notag
\end{gather}
The inequality in~\eqref{cond:contr-lem-comp-1} holds due to the bound $g(x,u) +\bfone/m\succ_x f(x,u)$. Since the system $\dot x = g(x, u)+\bfone/m$ is monotone, the inequality in~\eqref{cond:contr-lem-comp-2}  holds as well according to the remark after Proposition~\ref{prop:kamke}. This chain of inequalities contradicts~\eqref{cond:contr}, hence for all $t$ we have that
$\phi_m(t; x_2^m, u_2) \gg_x \phi_{f}(t; x_1, u_1)$. With $m\rightarrow+\infty$, by continuity of solutions we obtain $\phi_g(t, x_2, u_2) \succeq_x \phi_{f}(t; x_1, u_1)$, which completes the proof.
\end{pf}
\begin{proof-of}{Theorem~\ref{thm:comp-sys}} 
\textbf{A.} First we note that the assumption in~\eqref{eq:condss} implies that $s^0_g\preceq_x s^0_f\preceq_x s^0_r$. Indeed, take $x_0$ from the interior of the intersection of the sets $\cA(s^0_g)$, $\cA(s^0_f)$, $\cA(s^0_r)$. By Lemma~\ref{lem:comp-prin-control} for all $t$, we have $\phi_g(t, x_0, 0) \preceq_x \phi_f(t, x_0, 0) \preceq_x \phi_r(t, x_0, 0)$, and thus taking the limit $t\to\infty$ we get $s_g^0 \preceq_x s_f^0\preceq_x s^0_r$.\\
{\bf B.} Next we show that $g(x,u)\preceq_x f(x,u)$ for all $(x,u)\in\cD_M\times\cU$ implies that $\cS_{g}^-\supseteq\cS_{f}^-$. Let the set $\cV$ be such that $u=\mu h(\cdot,\tau)\in \cV$ if $(\mu,\tau)\in\cS^-_f$. \\
Due to $s_g^0 \preceq_x s_f^0$ and $g \preceq_x f$ on $\cD_M\times\cU_\infty$, by Lemma~\ref{lem:comp-prin-control}, we have that $s^0_g \preceq_x \phi_{g}(t; s^0_g, u) \preceq_x \phi_{f}(t; s^0_f, u )$, for all $u\in\cV$. Note that the first inequality is due to monotonicity of the system $\dot x = g(x,u)$. The flow $\phi_{f}(t; s^0_f, u)$ converges to $s^0_f$ with $t\rightarrow+\infty$. Therefore, there exists a time $T$ such that for all $t>T$ we have $s^0_g \preceq_x \phi_{g}(t; s^0_g, u ) \ll_x s^0_f+\varepsilon \bfone$ for some positive $\varepsilon$. Moreover, we can pick an $\varepsilon$ such that $s^0_f+\varepsilon \bfone$ lies in $\cA(s^0_g)$ (due to~\eqref{eq:condss}). Since the system $\dot x = g(x,u)$ is monotone, according to Lemma~\ref{lem:box-mon}, the flow $\phi_{g}(t; s^0_g, u)$ lies in $\cA(s^0_g)$. Hence, no $u$ in $\cV$ toggles the system $\dot x = g(x,u)$ either and we conclude that $\cS_{g}^-\supseteq\cS_{f}^-$. The proof that $\cS_{g}^-\supseteq\cS_{r}^-$ follows using the same arguments as above.\\
\textbf{C.} Finally, we show that $\cS_{f}^-\supseteq\cS_{r}^-$. Let the set $\cW$ be such that $u=\mu h(\cdot,\tau)\in\cW$ if  $(\mu,\tau)\in\cS_{r}^-$.\\
Due to $s_g^0 \preceq_x s_f^0 \preceq_x s_r^0$ and $g \preceq_x f \preceq_x r$ on $\cD_M\times\cU_\infty$, by Lemma~\ref{lem:comp-prin-control}, we have that 
\[\phi_{g}(t; s^0_g, u )\preceq_x 
\phi_{f}(t; s^0_f,u ) \preceq_x \phi_{r}(t; s^0_f, u), 
\]
for all $u\in\cW$. Now, monotonicity of $\dot{x}=g(x,u)$ implies that
$s^0_g \preceq_x \phi_{g}(t; s^0_g, u)$. Furthermore, there exists a $T$ such that $s^0_g \preceq_x  \phi_{f}(t; s^0_f, u) \preceq_x s^0_r+\varepsilon \bfone$ for all $t>T$, for all $u\in\cW$ and some small positive $\varepsilon$. This is due to the fact that $\phi_{r}(t; s^0_f, u)\rightarrow s^0_r$ with $t\rightarrow+\infty$. We can also choose an $\varepsilon$ such that $s^0_r+\varepsilon \bfone$ lies in $\cD_M$ due to~\eqref{eq:condss}. Hence, the flow of $\dot x = f(x,u)$ for all $u\in\cW$ belongs to the set $\{z | s^0_g \preceq_x  z \preceq_x s^0_r+\varepsilon \bfone\}$ for all $t>T$. \\
Now, assume there exists $u^c\in\cW$ that toggles the system $\dot x = f(x,u)$. This implies that the flow $\phi_{f}(t; s^0_f, u^c)$ converges to $s^1_f$ with $t\rightarrow\infty$. Therefore we have that $s^1_f$ belongs to the set $\{z | s^0_g \preceq_x  z \preceq_x s^0_r+\varepsilon \bfone\}$ for an arbitrarily small $\varepsilon$, and consequently $s^1_f\preceq_x s^0_r$. This contradicts the condition~\eqref{eq:condss2} in the premise of Theorem~\ref{thm:comp-sys}. Hence, no $u$ in $\cW$ toggles the system $\dot x = f(x,u)$ and $\cS_{f}^-\supseteq\cS_{r}^-$.
\end{proof-of}
\end{document}